\documentclass{article}

\usepackage{graphicx}%
\usepackage{multirow}%
\usepackage{amsmath,amssymb,amsfonts}%
\usepackage{amsthm}%
\usepackage{mathrsfs}%
\usepackage[title]{appendix}%
\usepackage{xcolor}%
\usepackage{textcomp}%
\usepackage{manyfoot}%
\usepackage{booktabs}%
\usepackage{algorithm}%
\usepackage{algorithmicx}%
\usepackage{algpseudocode}%
\usepackage{listings}%

\usepackage{pgfplotstable}
\usepackage{pgfplots}
\pgfplotsset{compat=1.16}

\usepackage{booktabs}
\usepackage{tikz}
\usetikzlibrary{patterns}
\usetikzlibrary{spy}
\usepackage{mathtools}
\usepackage{siunitx}

\usepackage{tabstackengine}
\stackMath

\usepackage{subcaption}
\usepackage{amsmath,amsfonts}
\usepackage{nicefrac}

\newcommand{\R}{\mathbb{R}}
\newcommand{\Ogrande}{\mathcal{O}}

\newcommand{\OSC}{\ensuremath{\mathcal{I}}}
\newcommand{\RBF}{\ensuremath{\mathsf{RBF}}}

\newcommand{\WENO}{\ensuremath{\mathsf{WENO}}}

\newcommand{\Pone}{\ensuremath{\mathsf{P1}}}
\newcommand{\CWENO}{\ensuremath{\mathsf{CWENO}}}
\newcommand{\CW}{\ensuremath{\mathsf{CW}}}

\newcommand{\stencil}{\mathcal{S}}
\newcommand{\rec}{\mathcal{R}}
\newcommand{\neighs}{\mathcal{N}}
\newcommand{\neighsFull}{\overline{\mathcal{N}}}
\newcommand{\dxmin}{\dx_{\text{min}}}
\newcommand{\doneSet}{\mathcal{D}}
\newcommand{\intd}{\ensuremath{\mathrm{d}}}
\newcommand\modulo[1]{\left\lvert#1\right\rvert}

\newcommand{\grid}{\mathcal{Q}}
\newcommand{\dt}{\mathrm{\Delta}t}
\newcommand{\dx}{\mathrm{\Delta}x}
\newcommand{\ddelta}{\mathrm{\Delta}}

\newcommand{\pcloud}{\mathcal{S}}

\newcommand{\pder}[2]{{#1}_{#2}}

\usepackage{algorithm}

\usepackage{graphicx}        %
\usepackage{multicol}        %
\usepackage{multirow}        %
\usepackage{array}
\usepackage{float}
\usepackage{hhline}
\usepackage[hang]{footmisc}
\newcolumntype{P}[1]{>{\centering\arraybackslash}p{#1}}

\usepackage{makecell}

\usepackage{url}

\theoremstyle{thmstyleone}%
\theoremstyle{thmstyletwo}%

\theoremstyle{thmstylethree}%
\newtheorem{definition}{Definition}%

\begin{document}
\title{Implicit reconstruction from point cloud: an adaptive level-set-based semi-Lagrangian method} %
\author{Silvia Preda%
 \thanks{Dipartimento di Scienza e Alta Tecnologia, Universit\`a dell'Insubria, Via Valleggio 11, 22100 Como (Italy)
   {\tt e-mail: silvia.preda@uninsubria.it} - ORCID 0009-0003-8405-2245}
 \and Matteo Semplice%
 \thanks{Dipartimento di Scienza e Alta Tecnologia, Universit\`a dell'Insubria, Via Valleggio 11, 22100 Como (Italy)
   {\tt e-mail: matteo.semplice@uninsubria.it} - ORCID 0000-0002-2398-0828}}
\maketitle

\begin{abstract}
We propose a level-set-based semi-Lagrangian method on graded adaptive Cartesian grids to address the problem of surface reconstruction from point clouds. The goal is to obtain an implicit, high-quality representation of real shapes that can subsequently serve as computational domain for partial differential equation models. The mathematical formulation is variational, incorporating a curvature constraint that minimizes the surface area while being weighted by the distance of the reconstructed surface from the input point cloud. Within the level set framework, this problem is reformulated as an advection–diffusion equation, which we solve using a semi-Lagrangian scheme coupled with a local high-order interpolator. Building on the features of the level set and semi-Lagrangian method, we use quadtree and octree data structures to represent the grid and generate a mesh with the finest resolution near the zero level set and the point cloud data. The complete surface reconstruction workflow is described, including localization and reinitialization techniques, as well as strategies to handle complex and evolving topologies. A broad set of numerical tests in two and three dimensions is presented to assess the effectiveness of the method.


\end{abstract} %

\section{Introduction}
\label{sec:intro}
Many real-world applications share the issue of acquiring and processing 3D digital models of non-synthetic objects. One of the common fundamental problem consists in reconstructing a smooth surface from a given set of unorganized points, called as \textit{point cloud}. This is typical, for instance, in fields like cultural heritage, where one deals with real artistic manufacts presenting complicated geometries and topologies, and whose shape is usually acquired, in the form of a point cloud, from the work of art itself, e.g. via 3D laser scanning or photogrammetry \cite{scan3d2020,Re:2011,ReEl:2006}. 

The importance of having a mathematical description of real shapes also lies in the problem of monitoring and predicting damage and degradation of monuments, for which some mathematical models, see for instance \cite{ADN:sulfation,CFN:08:swelling,CdFN:14:copper}, can be applied. These models require the solution of reaction-diffusion Partial Differential Equations (PDEs) on a computational domain having the exact shape of a work of art. This paper addresses the problem of providing a high fidelity description of a complex object to be later used as domain definition in PDE computations, in particular via ghost-cell methods \cite{OshFed:2004,GF:2002,CR2018,CoRu:24}. As a consequence, we seek a watertight representation, having some smoothness properties, which is suitable not only for static operations, but also for dynamic ones. Computations will involve high order methods and grids finer than the resolution of the point cloud, so that the reconstructed object can then be reliably represented on the extremely fine grids required by the PDE model.

In general, a surface can be represented by following two approaches: the explicit and the implicit one. Regarding the first, we mention mesh-based \cite{Am:delaunay2004,CaFr:surveyDelaunay} and parametric techniques, e.g. NURBS \cite{NURBS1,NURBS2}, which, despite their popularity in a wide variety of applications, are well-known to be difficult to handle when an evolution of the surface occurs. On the other hand, an implicit approach offers a better handling of topological flexibility, while having a very simple data structure that allows simple Boolean operations on the detected surface, too. 
Radial Basis Functions ($\RBF$) with global and local support have found many applications along this line \cite{RBF1,RBF2,RBF3} together with some least-squares-based methods \cite{leastSquare2}. Still, a signed distance function representing a surface can be obtained by using local estimators that associate an oriented plane to each point in the cloud \cite{hoppe1992}. A different and widespread approach for the processing of evolving surfaces, and point cloud data in particular, is constituted by the application of the Level Set Method (LSM) \cite{OshSet:1988}.   

Introduced in \cite{OshSet:1988}, the LSM has emerged as a powerful and versatile tool in a wide range of applications \cite{SSO:1994,Se:1999,OshFed:2004} including image processing and surface reconstruction \cite{Zhao:2000,LSApp:2004}. In its framework, an $n$-dimensional object $\Omega$ and its boundary $\Gamma$ are represented by a so called level set function $\phi:\R^n\rightarrow\R$ such that $\Omega = \{ \vec{x}\in\R^n : \phi(\vec{x})<0\}$ and $\Gamma$ is the zero isocontour of $\phi$, namely $\Gamma = \{ \vec{x}\in\R^n : \phi(\vec{x})=0\}$. Exploiting their representation via the scalar function $\phi$, $\Omega$ and $\Gamma$ can be then evolved in time accordingly to a PDE defined for $\phi$, driving the deformation of an initial surface. 

Along this line, a level-set-based semi-Lagrangian method with local interpolator has been presented in \cite{PreSe:25} for the problem of surface reconstruction from point cloud, while a first complete workflow from the dataset to PDE computations for the marble sulfation phenomenon can be found in \cite{cdss:mach19,CSS:monum}. 

The aforementioned method numerically solves the model introduced by Zhao et al. \cite{Zhao:2000}, which is based on the minimization of the energy functional
\begin{equation}\label{eq:energy}
E_p(\Gamma) = \Big( \int_{\Gamma} d^p(\vec{x})\, \intd s \Big) ^{1/p}, \quad 1 \leq p \leq \infty,  
\end{equation}
where $d(\vec{x}) = \min_{\vec{q} \in \pcloud}|\vec{x}-\vec{q} |$ is the distance function of a point $\vec{x}\in\R^n$ from the point cloud $\pcloud=\{\vec{q}_1,\ldots,\vec{q}_M\}$, and $\Gamma$ is a closed surface of co-dimension one in $\R^n$. The minimum of \eqref{eq:energy}, namely the final shape, is obtained by continuously deforming an initial surface $\Gamma^0$, following the gradient descent of \eqref{eq:energy}.

Within the LSM, the evolving surface $\Gamma(t)$ is represented implicitly using a level set function to capture the moving interface, leading to the level set formulation
\small
\begin{equation}\label{eq:levelset:pde}
    \begin{aligned}
        \pder{\phi}{t}(\vec{x},t) 
        &= \Bigg[  \frac{d(\vec{x})}{E_p(\phi)}   \Bigg]^{p-1} \Bigg( \nabla d(\vec{x})\cdot \nabla \phi(\vec{x},t) + \frac{\mu}{p} d(\vec{x}) \nabla\cdot\left( \frac{\nabla \phi(\vec{x},t)}{|\nabla \phi(\vec{x},t)|} \right)|\nabla \phi(\vec{x},t)| \Bigg),\\
        \phi(\vec{x},0) &= \phi^0(\vec{x}),
    \end{aligned}
\end{equation}
\normalsize
where the energy functional \eqref{eq:energy} is coherently rewritten as
\begin{equation} \label{eq:energy:ls}
    E_p(\phi) = \Bigg( \int_{\R^n} \lvert d(\vec{x}) \rvert ^p \delta (\phi) \lvert \nabla \phi \rvert \,\intd \vec{x} \Bigg)^{1/p},
\end{equation}
being $\delta$ the Dirac-delta function, and $\phi^0(\vec{x})$ a suitable initial data such that $\Gamma^0=\{\vec{x}\in\R^n: \phi^0(\vec{x})=0\}$. In \eqref{eq:levelset:pde}, the evolution of the surface towards the dataset $\pcloud$ is driven by the term $\nabla d(\vec{x})\cdot \nabla \phi(\vec{x},t)$, while the second term tempers the maximal curvature of $\Gamma$. The balance between these two is tuned by the parameter $\mu$, introduced in \cite{PreSe:25} for this purpose and to deal with noisy data, too. Finally, among all the possible level set representation of the surface $\Gamma$, we aim to work with the signed distance one, thus taking into account suitable reinitializazion procedure that guarantees our function $\phi$ to stay well-behaved, namely $\modulo{\nabla\phi}\approx 1$, during the evolution. For a wide description of the model and the parameters acting in \eqref{eq:levelset:pde} the reader can refer to \cite{Zhao:2000,Kosa2017,PreSe:25}. 

The numerical solution of \eqref{eq:levelset:pde} can be addressed in different ways, see for instance \cite{ZhOsFe:01,Kosa2017,CaFe:2017,He:2019}. In particular, in \cite{CaFe:2017}, a semi-Lagrangian (SL) scheme is adopted with the aim of overcoming the prohibitive time step constraint imposed by the diffusion term. In fact, first introduced in \cite{CoIsRe:52} for first-order systems of linear equations, SL schemes have gone trough an important extension with the main purpose of obtaining methods which are unconditionally stable with respect to the choice of the time step (see \cite{FaFe:14} for a comprehensive explanation), even in presence of parabolic terms, whose treatment has been studied in \cite{BF:14:SLdiffusion,BCCF:21}.

In the level set context, the possibility to overcome the parabolic-type CFL restriction is particularly tempting because it allows one to employ very fine grids for the computations, especially in the vicinity of the zero level set of $\phi$, also from an Adaptive Mesh Refinement (AMR) perspective. While the CFL condition will frustrate the adaptive approach when a general explicit scheme is employed, this will not be an issue in the SL framework where one is allowed to work at large Courant numbers, with a hyperbolic CFL restriction of type $\dt=\Ogrande(\dx)$. For more details, we refer the reader to \cite{St:99}, where a first application of SL schemes in the LSM context can be found.

Motivated by the reasons above, here we aim to extend the method presented in \cite{PreSe:25} to a fully adaptive framework, exploiting both the unconditional stability of the SL approach and the adaptivity criterion offered by the level set function itself. In fact, given a level set function $\phi$ associated to an interface $\Gamma$, constructing a tree-based grid refined around $\Gamma$ is quite simple since the values of $\modulo{\phi}$ naturally give a proper criterion for refinement, yielding a smaller grid size close to $\Gamma$ and a larger one far away from the interface. Moreover, we will couple this criterion with the proximity to the point cloud $\pcloud$, expressed by the distance function. As a result, the multiple iterations of the method presented in \cite{PreSe:25} are abandoned, in favour of a single computation in which the resolution of the computational grid increases along the zero level set and as we get closer to $\pcloud$.

On the other side, SL algorithms, compared to Eulerian ones, might not be easy to parallelize, especially on adaptive grids, since, depending of the CFL number, the feet of the characteristics may end up outside the halo region, in remote ranks, and may also significantly cluster, affecting load balancing. In order to reduce the computational effort and communications among processors, the SL scheme is coupled with two types of local reconstructions, a multilinear ($\Pone$) and a third-order accurate Central Weighted Essentially Non-Oscillatory ($\CWENO$) one, both based on a least-square approach due to the lack of structure in the mesh employed. 

First pioneering studies related to central reconstructions can be found in \cite{LPR:00:SIAMJSciComp}, where the authors suggested to blend, in a $\WENO$-like fashion, polynomials of different degrees, allowing to overcome some difficulties of non-existence, non-positivity and dependence on the reconstruction point of the $\WENO$ linear weights. The idea has been further developed into the so-called $\CWENO$ reconstruction and exploited in multiple spatial dimensions, also in the case of adaptive mesh refinement and non-uniform grids
\cite{BGFB:2020,ZhuQiu:2016,Baeza:19:CWENOglobalaverageweight,ZhouCai:08,ADER_CWENO}.
The technique has also been exploited in finite difference schemes for Hamilton-Jacobi (HJ) equations on Cartesian meshes via dimensional splitting \cite{ZhuQiu:2017:HJ,ZhengShuQiu:2019:HJ}, as well as on general meshes \cite{ZhuQiu:2020:triHJ}.
General results for establishing the convergence order of a $\CWENO$ reconstruction have been presented in \cite{CPSV:cweno,CSV19:cwenoz,SV:2020:CWAO}.

In order to design the complete numerical method and achieve the final reconstruction via a signed distance function, some other techniques will be employed, including localization procedure \cite{PENG1999}, reinitialization \cite{SSO:1994,OshSet:1988,AMR:Saye:14} and distance function computation \cite{Zhao2005AFS}. 

The paper is organized as follows. The SL scheme for the numerical approximation of \eqref{eq:levelset:pde} is described in \S\ref{sec:numericalscheme}. All the auxiliary schemes needed for the level set computation are collected in \S\ref{sec:auxiliaryschemes}. In \S\ref{sec:numericaltests}, the whole algorithm is summarized and various numerical tests in two and three dimensions are presented. Finally, \S\ref{sec:conclusions} collects some conclusions and future perspectives.

\section{Numerical scheme}
\label{sec:numericalscheme}

We consider a background quadtree (resp. octree) mesh having a single tree-based structure which is defined on a cubic domain $\Omega'\subset\R^2$ (resp. $\R^3$) containing the point cloud $\pcloud$ and $\Omega(t)$, at each time $t>0$. In what follows, the set of quadrants composing $\Omega'$ will be denoted by $\grid$ and, for $j\in\grid$, $\phi_j^n$ indicates the approximate value of $\phi$ evaluated at the centre $\vec{x}_j$ of the quadrant $j$ at time $t^n$.
Moreover, for each quadrant $j$, we denote with $\dx_j$ its edge length and, once the maximum level of refinement $L$ is fixed, we define as
\begin{equation}\label{eq:AMRdxMin}
    \dx_{\text{min}} = \frac{\ddelta\Omega'}{2^L}
\end{equation}
the minimum edge length of the quadrants composing the tree, where $\ddelta\Omega'$ is the edge length of the cubic domain $\Omega'$.

The discretized values of $\phi$, namely the numerical evolution of \eqref{eq:levelset:pde}, are computed following the SL approach presented in \cite{CaFe:2017}, coupled with a suitable local $\Pone$ or $\CWENO$ reconstruction to handle the quadtree grid, reduce the computational costs and exploit the parallelization. The SL scheme will be briefly described in the next \S~\ref{ssec:semilagrangian}, for the two dimensional case, while the reader can refer to \cite{CaFe:2017,PreSe:25} for a complete description in three dimensions. 
The complementary techniques and the complete algorithm will be described later on, in \S~\ref{sec:auxiliaryschemes} and \S~\ref{sec:numericaltests}, respectively.

\subsection{Semi-Lagrangian scheme}\label{ssec:semilagrangian}

In two dimensions, we compute the update of $\phi^n_j$ as
\begin{equation}\label{eq:SL2d}
    \begin{aligned}
       \phi^{n+1}_j
        \, &=\; \frac12 \sum_{i=1}^2\rec[\Phi^n]\left( \vec{x}^*_{j,i} \right),
        \\
        \vec{x}^*_{j,i} &= \vec{x}_j + C^n_j \dt \nabla d(\vec{x}_j) + \sqrt{\frac{2 C^n_j \, \mu \, d(\vec{x}_j) \dt}{p}}\,\vec{\sigma}^n_j \lambda_i,
    \end{aligned}
\end{equation}
where $\lambda_i$ ranges over $\{ -1, +1\}$ and $C^n_j$ is the scale factor $\big[\frac{d(\vec{x_j})}{E_p(\phi^n)}\big]^{p-1}$. The operator $\rec[\Phi^n](\vec{x})$ denotes a reconstruction at point $\vec{x}$ of the data $\Phi^n=\{\phi^n_j: j\in\grid\}$, which will be specified later in \S~\ref{ssec:P1_LS} and \S~\ref{ssec:CW_LS}. The unit vector $\vec{\sigma}^n_j$ is given by
\begin{equation}\label{sigma2d}
    \vec{\sigma}^n_j = \frac{1}{\lvert\nabla\phi^n_j\rvert} 
    \begin{pmatrix} 
        &\partial_y\phi^n_j \\ -&\partial_x\phi^n_j
    \end{pmatrix},
\end{equation}
is tangent to the level sets of $\phi$ and thus orthogonal to its gradient.

One can notice that, in equation \eqref{eq:SL2d}, the evaluation points of the reconstruction operator are obtained by computing the foot of the characteristic pertaining to the advection term in \eqref{eq:levelset:pde}, and adding a further displacement of size $\sqrt{\dt}$ that generates a diffusion along the tangent space of the level sets, thus discretizing the curvature term in \eqref{eq:levelset:pde}, as described in \cite{FaFe:2003,CaFaFe:2010}.

Because of the factor $1/\lvert\nabla\phi^n_j\rvert$ inside the definition of $\vec{\sigma}^n_j$ in \eqref{sigma2d}, when $|\nabla \phi^n_j|< D\dt^\alpha$, the scheme is replaced by
\begin{equation}\label{eq:SL:2}
\phi^{n+1}_j = \frac{1}{\lvert \neighs_j \vert} \sum_{i\in\neighs_j} \rec[\Phi^n](\vec{x}_i), 
\end{equation}
where $\neighs_j$ is the set of the first neighbour indexes of $j\in\grid$ and $\lvert \neighs_j \vert$ represents its cardinality.
In all computations presented in this work, we set $D=10^{-3}$ and $\alpha=1$.

Analogous computations can be done in order to obtain the 3D version of the scheme presented above. For its complete description we refer the reader to \cite{CaFe:2017,PreSe:25}.

\paragraph{Remark}
In order to apply the SL scheme \eqref{eq:SL2d}, one has to quickly find the quadrant to which the points $\vec{x}^*_{j,i}$ belong. While on Cartesian grids localizing them is a trivial task, on quadtree this is no longer an easy procedure as it involves browsing the tree structure from the root until one finds the leaf containing the point $\vec{x}^*_{j,i}$.\\

In order to fully detail the scheme \eqref{eq:SL2d}, a suitable reconstruction operator has to be defined. Considering a space of polynomials $\mathbb{P}$ with a basis $\mathcal{B}$, the key point of interpolation on quadtree is that, since each quadrant $j\in\grid$ may have a different number of neighbours $\modulo{\neighs_j}$, the linear system arising for each quadrant $j$ may have a number of rows other than $\modulo{\mathcal{B}}$. Thus, in our scheme we will resort to a $\Pone$ constrained least-squares reconstruction $\rec_{P1}[\Phi]$ and to a third order $\CWENO$ one $\rec_{CW}[\Phi]$, based on the same constrained least-squares approach.

\begin{figure}[htpb]
\centering
\begin{minipage}{0.45\linewidth}
    \centering
    \includegraphics[width=1\linewidth]{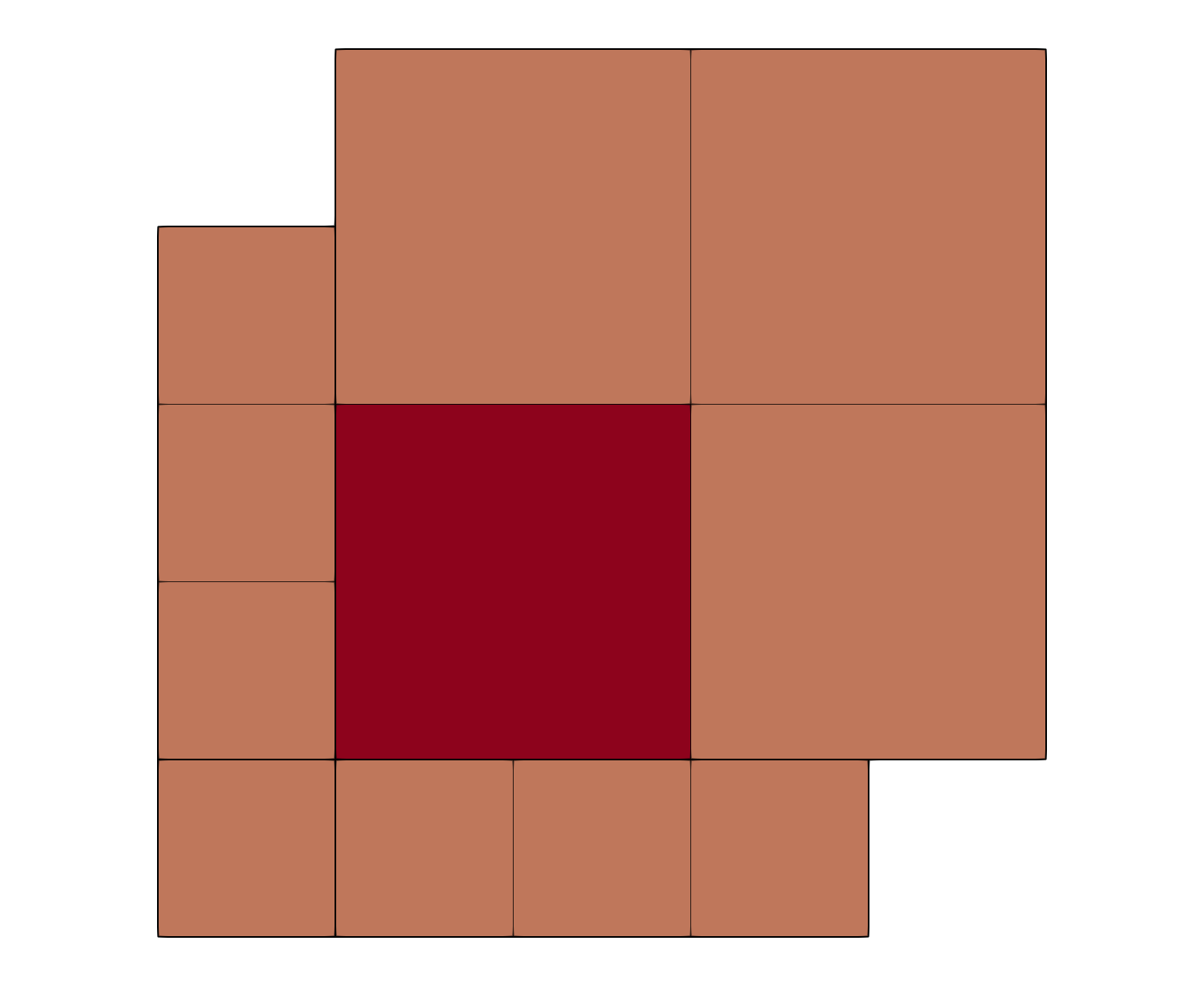}
\end{minipage}
\begin{minipage}{0.45\linewidth}
    \centering
    \includegraphics[width=1\linewidth]{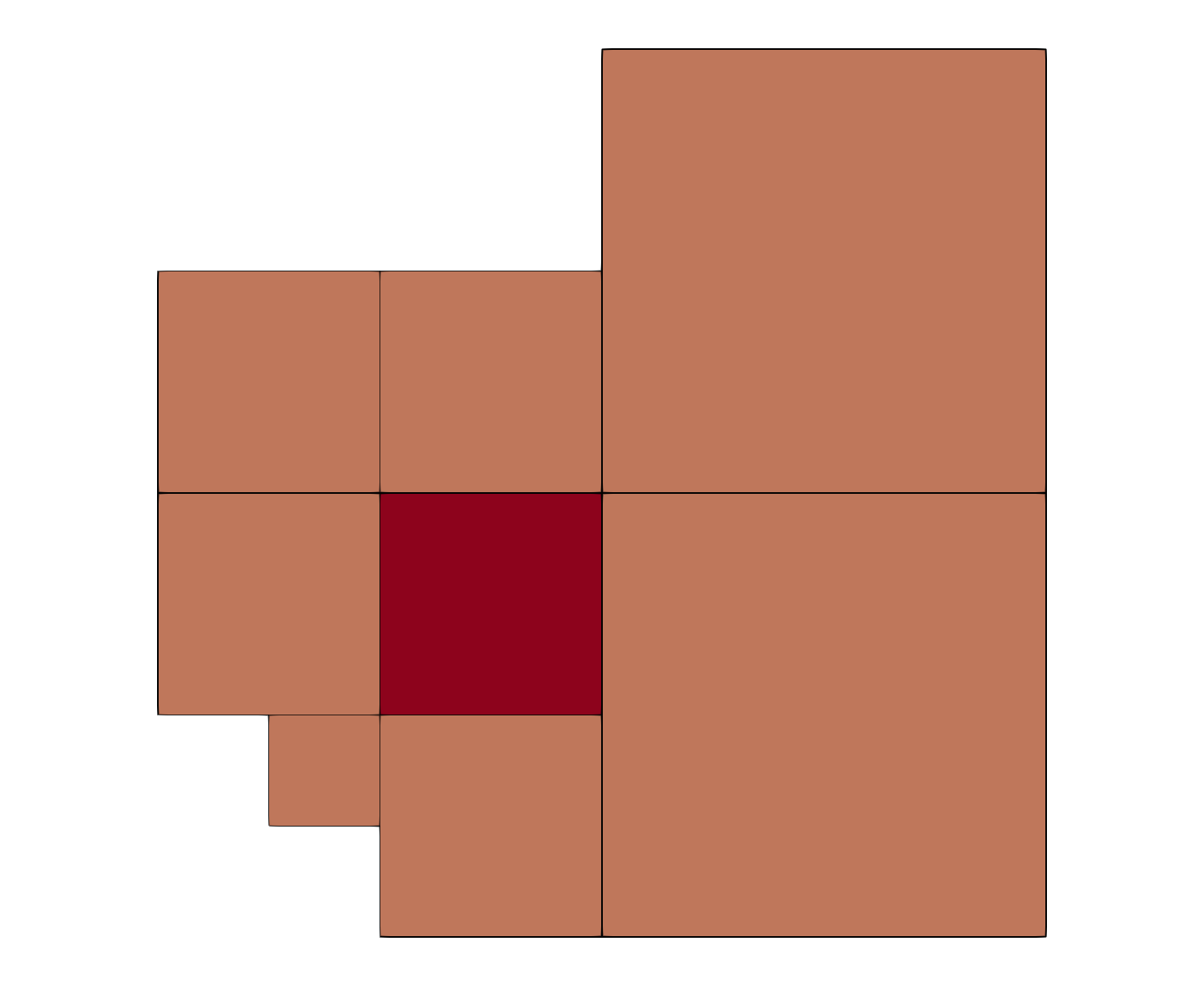}
\end{minipage}
\caption{Example of quadrants with different neighbours. The main quadrants and their neighbours are depicted in red and orange, respectively.} 
\label{fig:Q1stencils}	
\end{figure}

\subsection{$\Pone$ constrained least-squares reconstruction} 
\label{ssec:P1_LS}

In 2D, on a quadrant $j$ of size $\dx$, we consider the local basis 
$\mathcal{B}_j = \{\varphi_{j_s}\}_{s=0}^2 = \{ 1,\hat{x},\hat{y}\}$, where $\hat{x}=(x-x_j)/\dx$ and $\hat{y}=(y-y_j)/\dx$. We seek for the vector of coefficients $\{c_{j_s}\}_{s=0}^2$ by imposing the constraint $c_{j_0}=\phi_j$ and the interpolation conditions
\begin{equation}\label{eq:P1_LS_condition}
\sum\limits_{s=1}^2 c_{j_s}\varphi_{j_s}(\vec{x}_i) = \phi_i-\phi_j, \quad i\in\neighs_j,
\end{equation}
such that $\rec_{P1}[\Phi](\vec{x}_j) = \phi_j$ holds true.

With this choice one can immediately note that, considering both edge and corner neighbours (see Figure~\ref{fig:Q1stencils}), the case $\modulo{\neighs_j}<\modulo{\mathcal{B}_j}-1$ never occurs. Instead, the system for the $c_{j_s}$ is overdetermined and is solved using a least-squares approach. 

In an analogous way, we can compute the $\Pone$ reconstruction for a quadrant $j$ in the 3D case, choosing $\mathcal{B}_j = \{ 1,\hat{x},\hat{y},\hat{z}\}$.\\

Despite its low computational cost, we expect that the $\Pone$ operator might lead to less accurate reconstructed surfaces than alternative high-order ones. However, high-degree polynomials might introduce oscillations that severely affect the stability of the scheme. Essentially non-oscillatory techniques are apt for this purpose, producing high-order non-oscillatory operators by introducing a non-linear dependency from the data. In particular, we resort to a $\CWENO$ operator which blends and selects first degree polynomials, having smaller and directionally biased stencils, with the optimal second degree one built on the entire set of neighbours of a quadrant and the quadrant itself, in order to get high-order accuracy. Due to the lack of structure in the adaptive mesh employed, both optimal and lateral polynomials will be obtained with the same constrained least-squares approach presented for the $\Pone$ case.

\subsection{$\CWENO$ constrained least-squares reconstruction} 
\label{ssec:CW_LS}

Let $j\in\grid$ be the quadrant in which we want to compute the reconstruction. We recall now the general definition of the $\CWENO$ operator, along the lines of \cite{CraSeVi:19:cwenoz}. 

\begin{definition}
Given a stencil $\stencil_{\text{opt}}$, including $j$, let $P_{\text{opt}}\in \mathbb{P}^G_n$ be the optimal polynomial of degree $G$, associated to $\stencil_{\text{opt}}$. Further, let $P_1,P_2,...,P_m$ be a set of $m\geq 1$ polynomials of degree $g$ with $g<G$, associated to substencil $\stencil_k$ such that $\{j\}\subset \stencil_k \subset \stencil_{\text{opt}}$ $\forall k=1,...,m$. Let also $\left\{ d_k \right\}_{k=0}^m$ be a set of strictly positive real coefficients such that $\sum_{k=0}^m d_k = 1$.

The $\CWENO$ operator computes a reconstruction polynomial 
\begin{equation}
    \rec_{CW} =\CWENO(P_{\text{opt}},P_1,...,P_m) \in \mathbb{P}^G_n,
\end{equation}
as follows:
\begin{enumerate}
    \item introduce the polynomial $P_0\in \mathbb{P}^G_n$ defined as
    \begin{equation}
        P_0(x) = \frac{1}{d_0}\left( P_{\text{opt}}(x) - \sum_{k=1}^m d_k P_k(x) \right);
    \end{equation}
    \item compute suitable regularity indicators
    \begin{equation}
        \OSC_0 = \OSC[P_{\text{opt}}], \quad \OSC_k = \OSC[P_k], \quad k\geq 1;
    \end{equation}
    \item compute the nonlinear coefficients $\left\{ \omega^Z_k \right\}_{k=0}^m$ as
        \begin{equation}\label{eq:CWENOomega}
        \alpha_k^Z = \frac{d_k}{(\OSC_k + \epsilon)^l}, \quad \omega_k^{CW} = \frac{\alpha_k^Z}{\sum_{i=0}^m \alpha_i^Z}
        \end{equation}
    where $\epsilon$ is a small positive quantity and $l\geq 1$;\\
    \item define the reconstruction polynomial as
    \begin{equation}\label{eq:CWENOrec}
        \rec_{CW}(x) = \sum_{k=0}^m \omega_k^{CW} P_k(x) \in \mathbb{P}^G_n.
    \end{equation}
\end{enumerate}
\end{definition}

Compared to the traditional $\WENO$ operator, the reconstruction polynomial defined in \eqref{eq:CWENOrec} can be evaluated at any point in the quadrant $j$ at a very low computational cost. This has been proven to be an advantage in SL schemes (see \cite{CaFePreSe:25}) that require to evaluate the reconstruction in different points laying in the same quadrant, which is, indeed, the case of our scheme \eqref{eq:SL2d}. 

The accuracy and non-oscillatory properties of the $\CWENO$ scheme is guaranteed by the dependence of its non-linear weights \eqref{eq:CWENOomega} on the regularity indicators $\OSC_k$. On smooth data, the regularity indicators of all polynomials will be close to each others, so that $\rec\approx P_{\text{opt}}$ and the reconstruction reaches the optimal order of accuracy $G+1$. 
On the other hand, when a discontinuity is present in $\stencil_{\text{opt}}$, both $\OSC_0\asymp1$ and at least one $\OSC_{\hat{k}} \asymp 1$ for some $\hat{k}\in \{ 1,...,m \}$, where with $\asymp1$ we mean ``of the same order of $1$, and bounded away from zero''. In this case, the reconstruction is given by a linear combination of all lateral polynomials that are not affected by the discontinuity and its accuracy reduces to $g+1$.

\subsubsection{Two spatial dimensions}

Let $j\in\grid$ be the quadrant of the grid containing the reconstruction point $x$ and let $\neighsFull_j = \neighs_j\cup\{j\}$ be the stencil for the reconstruction, where $\neighs_j$ is the set of neighbouring quadrants of $j$.
The setup of the stencils for optimal and lateral polynomials is the following: the optimal polynomial is the one associated to $\neighsFull_j$, while the substencils for the south-west (sw), south-east (se), north-west (nw) and north-east (ne) polynomials are composed by $j$ itself and corner or edge neighbours in the corresponding directions. 

The reconstruction operator compute $\rec_{CW}\in \mathbb{P}^2_2$ as
\small
\begin{equation}
    \rec_{CW} = \CWENO(P_{opt}^{(2)}; P^{(1)}_{\text{sw}}, P^{(1)}_{\text{se}}, P^{(1)}_{\text{nw}}, P^{(1)}_{\text{ne}} ),
\end{equation}
\normalsize
with the optimal and lateral polynomials set to be second and first degree polynomials, respectively, built on the local basis
\begin{equation}\label{eq:basisCWENO2d}
    \mathcal{B}_2 = \{ 1,\hat{x},\hat{y},\hat{x}^2,\hat{x}\hat{y},\hat{y}^2\}, \quad
    \mathcal{B}_1 = \{ 1,\hat{x},\hat{y}\}.
\end{equation}
For each polynomial we impose the constraint $P(\vec{x}_j)=\phi_j$, such that $c_{j_0}=\phi_j$, and compute the remaining coefficients $c_{j_k}$ similarly to \eqref{eq:P1_LS_condition}, following a least-squares approach. Since we impose the constraint on each polynomial, we can guarantee that the condition $\rec_{CW}[\Phi](\vec{x}_j)=\phi_j$ holds true. 

Once the coefficients are obtained, the oscillation indicators are defined as
\small
\begin{equation}\label{eq:IND:2d}
\OSC[P] 
= 
\sum_{|\vec{\alpha}|\geq1} \dx_j^{2|\vec{\alpha}|-2} 
\int_{Q_j}
\left[
  \frac{\partial P}{\partial x^{\alpha_1}\partial y^{\alpha_2}}
\right]^2
\intd x\intd y,
\end{equation}
\normalsize
according to the classical ones of \cite{JiangShu:96,HuShu:99}. This choice might be different to the usual one for HJ equation, whose solutions can be at worst continuous with kinks, see \cite{JiPe:00:HJind,FaPaTo:20:HJind2d}, but here is motivated by the degrees of the polynomials employed in the reconstruction. In fact, since our reconstruction employs lateral polynomials of at most first degree, their HJ-type indicators would be all identically equal to zero, providing no information about their regularity and thus the ones to be discarded in case of oscillations, consequently affecting the non-oscillatory properties of the entire reconstruction procedure. 

Denoting with $\vec{c}$ the vector of coefficients of a generic polynomial in two variables of degree up to 2, along the basis $\mathcal{B}_2$ in \eqref{eq:basisCWENO2d}, its oscillation indicator can be expressed as
\begin{equation}\label{eq:ind2dCoeff}
\OSC[P] = \vec{c}^{\,\,\small T} M \vec{c},
\end{equation}
with
\begin{equation}
\setstackgap{L}{1.1\baselineskip}
\fixTABwidth{T}
M = \parenMatrixstack{
0 & 0 & 0 & 0 & 0 & 0 \\
0 & 1 & 0 & 0 & 0 & 0 \\
0 & 0 & 1 & 0 & 0 & 0 \\
0 & 0 & 0 & \nicefrac{13}{3} & 0 & 0 \\
0 & 0 & 0 & 0 & \nicefrac{7}{6} & 0 \\
0 & 0 & 0 & 0 & 0 & \nicefrac{13}{3} 
}.
\end{equation}

The $\CWENO$ reconstruction applied in the numerical tests of this paper is then defined by choosing linear coefficients $d_{\text{sw}}=d_{\text{se}}=d_{\text{nw}}=d_{\text{ne}}=\nicefrac{1}{8}$ and $d_0=\nicefrac{3}{4}$, $l=2$ and $\epsilon = \dx^2$ in the above construction.

\subsubsection{Three spatial dimensions}

The construction for the $\CWENO$ operator in three dimensions strictly follows the guidelines seen above for the two-dimensional case.

The reconstruction operator blends the second degree optimal polynomial associated to the stencil $\neighsFull_j$ with 8 first degree polynomials corresponding to the biased ones. To collect these stencils we consider west, east, north, south, backward and forward directions. The result of the reconstruction operator is $\rec_{CW}\in\mathbb{P}^2_3$, built on the basis 
\begin{equation}\label{eq:basisCWENO3d}
\begin{aligned}
    \mathcal{B}_2 = \{ 1,\hat{x},\hat{y},\hat{z},\hat{x}^2,\hat{x}\hat{y},\hat{y}^2,\hat{z}^2,\hat{x}\hat{z},\hat{y}\hat{z}\}, \\
    \mathcal{B}_1 = \{ 1,\hat{x},\hat{y},\hat{z}\}.
\end{aligned}
\end{equation}

Due to the dimension, the oscillation indicators need to be properly rescaled and are thus defined as
\small
\begin{equation}\label{eq:IND:3d}
\OSC[P] 
= 
\sum_{|\vec{\alpha}|\geq1} \dx_j^{2|\vec{\alpha}|-3} 
\int_{Q_j}
\left[
  \frac{\partial P}{\partial x^{\alpha_1}\partial y^{\alpha_2}}
\right]^2
\intd x\intd y\intd z,
\end{equation}
\normalsize
or in a matrix form as
\begin{equation}\label{eq:ind3dCoeff}
\OSC[P] = \vec{c}^{\,\,\small T} M \vec{c},
\end{equation}
with
\begin{equation}
\setstackgap{L}{1.1\baselineskip}
\fixTABwidth{T}
M = \parenMatrixstack{
0 & 0 & 0 & 0 & 0 & 0 & 0 & 0 & 0 \\
0 & 1 & 0 & 0 & 0 & 0 & 0 & 0 & 0 \\
0 & 0 & 1 & 0 & 0 & 0 & 0 & 0 & 0 \\
0 & 0 & 0 & \nicefrac{13}{3} & 0 & 0 & 0 & 0 & 0 \\
0 & 0 & 0 & 0 & \nicefrac{7}{6} & 0 & 0 & 0 & 0 \\
0 & 0 & 0 & 0 & 0 & \nicefrac{13}{3} & 0 & 0 & 0 \\ 
0 & 0 & 0 & 0 & 0 & 0 & \nicefrac{13}{3} & 0 & 0 \\
0 & 0 & 0 & 0 & 0 & 0 & 0 & \nicefrac{7}{6} & 0 \\
0 & 0 & 0 & 0 & 0 & 0 & 0 & 0 & \nicefrac{7}{6} 
}.
\end{equation}

\paragraph{Remark}

The scaling factor in \eqref{eq:IND:2d} and \eqref{eq:IND:3d} is motivated by the expected behaviour of the indicators. To clarify the idea, let us consider the simpler case of a generic polynomial $P_k$ of degree $r-1$, in one dimension, having a discontinuity in its stencil $S_k$. In this case, the indicator would read
\small
\begin{equation}\label{eq:IND:1d}
\OSC[P_k] 
= 
\sum_{\alpha\geq1} \dx_j^{2\alpha-1} 
\int_{Q_j}
\left[
  \frac{\partial P}{\partial x^{\alpha}}
\right]^2
\intd x.
\end{equation}
\normalsize
The $\alpha$-derivative of $P_k$ %
is a polynomial of degree $r-\alpha-1$, whose coefficients of degree $m$ %
are given by the Newton divided differences with $m+\alpha+1$ points, which scale as $\Ogrande(\dx^{-m-\alpha})$, due to the discontinuity. Computing the integral in \eqref{eq:IND:1d} contributes another $\Delta x$ factor and the indicator is thus a sum of terms, all scaling at the same rate $\Ogrande(\dx^{1-2\alpha})$. It is thus clear that the exponent $2\alpha-1$ in the scaling factor yields $\OSC_k\asymp1$ in this non-regular case; on the other hand, in case of smoothness, the Newton divided differences would be $\Ogrande(1)$ and one gets $\OSC_k\to 0$, as expected.

With similar arguments one can deduce the correct scaling factor in the definition of the indicator in the general multi dimensional case.

\section{Computation of the level set evolution}
\label{sec:auxiliaryschemes}

The overall surface reconstruction procedure requires different schemes to be put together. The SL scheme presented in the previous section constitutes the main ingredient of the method; all the other auxiliary parts will be detailed later on in this section. A summary of the complete algorithm will be given at the beginning of \S~\ref{sec:numericaltests}.

\subsection{Localization}
\label{ssec:localization}

A typical approach when dealing with LSM is to localize the evolution only in a narrow band around the zero level set of $\phi$. Along the lines of \cite{PENG1999}, the definition of the narrow band, and related computations, is based on the two constants
\begin{equation}\label{eq:beta}
    \beta = 2\lambda\dxmin,
\end{equation}
and
\begin{equation}\label{eq:gamma}
    \gamma = 4\lambda\dxmin,
\end{equation}
where $\dxmin$ is the minimum edge length of the refined grid, as defined by \eqref{eq:AMRdxMin}, and $\lambda = \dt/\dxmin$ is related to the choice of the time step $\dt=\Ogrande(\dxmin)$.

Since it is not important to update $\phi$ far away from zero, at each time step we choose a subgrid
\begin{equation}\label{eq:subgrid}
\widetilde\grid=\{j\in\grid: \modulo{\phi^n_j}<\gamma\}\subset\grid,     
\end{equation}
according to \cite{PENG1999} and to the stencils for the spatial reconstructions employed in this work, and update $\phi^{n+1}$ only therein. Outside $\widetilde\grid$, we simply cut our level set function as
\begin{equation}\label{eq:cut}
    \phi_j = 
    \begin{cases}
        \gamma \quad &\text{if} \quad \phi_j>\gamma,\\
        \phi_j \quad &\text{if} \quad \modulo{\phi_j}\leq\gamma,\\
        -\gamma \quad &\text{if} \quad \phi_j<-\gamma.
    \end{cases} 
\end{equation}

Analogously to \cite{PENG1999}, to prevent numerical oscillations at the boundary of $\widetilde\grid$, we update the solution involving the additional cut-off function \begin{equation}\label{eq:cutoff}
    c(\phi)=
    \begin{cases}
        1 \hfill &\text{if} \; |\phi|\leq\beta,\\
        \frac{(|\phi|-\gamma)^2(2|\phi|+\gamma-3\beta)}{(\gamma-\beta)^3} &\text{if} \; \beta<|\phi|\leq\gamma,\\
        0 \hfill &\text{if} \; |\phi|>\gamma,
    \end{cases}
\end{equation}
considering the modified equation
\small
\begin{equation}\label{eq:ZhaoPDECutOff}
\begin{aligned}
    \pder{\phi}{t}(\vec{x},t) &= c(\phi) \Bigg[  \frac{d(\vec{x})}{E_p(\phi)}   \Bigg]^{p-1} \Bigg( \nabla d(\vec{x})\cdot \nabla \phi(\vec{x},t) \\&\quad + \frac{\mu}{p}d(\vec{x}) \nabla\cdot\left( \frac{\nabla \phi(\vec{x},t)}{|\nabla \phi(\vec{x},t)|} \right)|\nabla \phi(\vec{x},t)| 
      \Bigg)
\end{aligned}
\end{equation}
\normalsize
and the proper modifications to the semi-Lagrangian scheme (see \cite{PreSe:25}). 

Finally, even if we consider the modified equation \eqref{eq:ZhaoPDECutOff}, keeping the level set function close to the signed distance one during its evolution constitutes an important issue and a reinitialization procedure will be required. This reinitialization procedure must be performed on a larger narrow band $\overline\grid$ obtained by considering $\widetilde\grid$ and a number of neighbouring layers proportional to $\lambda$. The fact that  $\overline\grid \supset \widetilde\grid$ is essential if one starts the algorithm with an initial data so far from the point cloud that the initial computational band \eqref{eq:subgrid} does not contain $\pcloud$. In such a case, without the enlarged reinitialization band, the evolution of $\phi$ would remain confined to the first computational band $\widetilde\grid^0$, while with this choice, the successive bands $\widetilde\grid^n$ will be able to move contextually with the zero level set of $\phi^n$. 

\subsection{Adaptivity}
\label{ssec:adaptivity}

We consider an adaptive framework in which the computational domain is subdivided into quadrants (resp. octants). Each of them can be refined by splitting into $4$ (resp. $8$) equal pieces, each group of them can be coarsened and replaced by the quadrant (octant) they were obtained from. Each quadrant has a level $l\geq 0$, which dictates its edge length as being $1/2^l$ that of the computational domain. The minimum and maximum level control the size of the largest and smallest cell in the grid. Typically one also enforces that the levels of two adjacent quadrants do not differ by more than $1$. The grids can be represented by a quadtree (resp. octree), in which any refined quadrant is a node with $4$ (resp. $8$) children and the active computational quadrants are the leaves of the tree.

In what follows, unless otherwise specified, we will use the terms cell, quadrant, octant and leaf, to refer to the elements of the computational grid.

In the AMR context the typical adaptivity procedure follows the steps below:
\begin{itemize}
    \item first of all, the leaves are marked for refinement, coarsening or to be left unchanged, depending on the chosen criteria;
    \item the refinement and coarsening is then applied to each marked leaf in a recursive way and consistently with the minimum and maximum level of refinement set for the mesh;
    \item if the grid is graded, once it has been roughly adapted, a 2:1 balancing is performed in order to guarantee that the level difference between a quadrant and each of its neighbours is at most 1;
    \item finally, as far as parallelization is concerned, load distribution is operated between processes by an equal division of the new array of leaves.
\end{itemize}

In our level set application solid adaptivity criteria need to be defined to ensure the grid to be properly refined in the vicinity of the zero level set, i.e. the surface interface $\Gamma$, and in the proximity of the cloud data $\pcloud$. A progressive decreasing resolution is expected when moving further away, in particular where our level set function $\phi$ is flattened to a prescribed value $\pm\gamma$, according to the localization of the method, described in \S~\ref{ssec:localization}.

\subsubsection{Refinement}

Let us indicate with $h_{\pcloud}$ a suitable estimation of the resolution of the point cloud $\pcloud$, i.e. the average distance between two closest points in $\pcloud$, and let us recall that $l$ and $L$ denote the current level of a quadrant and the maximum level of refinement, respectively. In order to refine the computational grid we proceed as follows:
\begin{itemize}
    \item a quadrant $j\in\grid$ is refined up to level $L$ if
\begin{equation}
    \modulo{\phi_j}<\gamma \quad \text{and} \quad d<2h_{\pcloud};
\end{equation}
    \item a quadrant $j\in\grid$ is refined up to level $L-1$ if
\begin{equation}
    \modulo{\phi_j}<\gamma \quad \text{and} \quad 2h_{\pcloud}\leq d<4h_{\pcloud};
\end{equation}
    \item a quadrant $j\in\grid$ is refined up to level $L-2$ if
\begin{equation}
    \modulo{\phi_j}<\gamma;
\end{equation}
\end{itemize}
where $\gamma$ is defined as in \S~\ref{ssec:localization}. 

Once the quadrant $j$ is marked, the procedure for the refinement is trivial thanks to the local reconstruction $\rec_j$ defined on the quadrant itself. The quadrant $j$ of centre $\vec{x}_j$ originates $4$ children in 2D, $8$ in 3D, by dividing by two each of its edges. The new quadrants $j_r$ have centres $\vec{x}_{j_r}$ and their corresponding level set values are $\phi_{j_r}=\rec_j[\Phi](\vec{x}_{j_r})$.

\subsubsection{Coarsening}

On the other side, a set of 4 in 2D (resp. 8 in 3D) children $C\subset\grid$ is coarsened if 
\begin{equation}\label{eq:coarsCriterion}
    \modulo{ \{ c\in C | \modulo{\phi_{c}}<\gamma \} } = 0,
\end{equation}
where $\gamma$ is again defined as in \S~\ref{ssec:localization}.  

The coarsening procedure consists in gluing together quadrants of the same level originated from the same branch: the new centre $\vec{x}_j$ and the corresponding level set value is obtained by averages from the outgoing quadrants of centres $\{\vec{x}_c\}_{c\in C}$, namely
\begin{equation}\label{eq:coarsProc}
    \vec{x}_j = \frac{1}{\modulo{C}}\sum\limits_{c\in C} \vec{x}_c
    \quad \text{and} \quad
    \phi_j = \frac{1}{\modulo{C}}\sum\limits_{c\in C} \phi_c.
\end{equation}

\subsection{Distance function computation}
\label{ssec:distance}

The evolution of the level set function $\phi$ requires, as a fundamental step, to compute the distance function $d$ from the point cloud $\pcloud$, which we recall to be given by $d(\vec{x}) = \min_{\vec{q} \in \pcloud}|\vec{x}-\vec{q} |$. The distance function appears itself in the governing equation \eqref{eq:levelset:pde}, in the definition of the energy functional, and its gradient constitutes the velocity field driving the evolution of the level set towards the data. Clearly, one is not interested in computing the exact values of $d$ on the whole grid, but just in the vicinity of the cloud.

Thus, in a first step, $d_j$ can be computed exactly in the quadrants that contain at least a point of $\pcloud$ as the minimum distance of the centre $\vec{x}_j$ from the cloud points located in the quadrant $j$. Then, this information can be spread across the whole grid, all over the processors, by solving the Eikonal equation
\begin{equation}\label{eq:distEikonal}
    \modulo{\nabla d(\vec{x})} = 1, \quad d(\vec{q}\in\pcloud)=0.
\end{equation}
The most popular algorithms for solving \eqref{eq:distEikonal} are the Fast Marching Method (FMM) \cite{FMM:Tsi:95,FMM:Se:96,FMM:Cho:01} and the Fast Sweeping Method (FSM) \cite{Zhao:2000,FSM:Tsai:03}. However, in our framework, these methods should be applied carefully since they might suffer both from the adaptive discretization of the grid and the parallelization of the scheme, due to the causality across processes.
Some of the earliest attempts in parallelizing the FMM were proposed in \cite{AMR:He:03,AMR:Tu:08}, where one could notice how the number of iterations needed to get convergence greatly depends on the complexity of the interface and on the parallel partitioning and, in general, fewer iterations are required if the domains are aligned with the normals to the interface. In \cite{AMR:Zhao:07} a parallel FSM method was presented for the first time, while a hybrid FMM-FSM was presented in \cite{AMR:FMM-FSM:15}. A parallel Fast Iterative Method (FIM) has been also proposed in \cite{AMR:FIM:08}. The case of unstructured meshes has been investigated in \cite{FSM:tri:Qian:07,FSM:ustrTri:24}.

Due to the complexity of the aforementioned methods on AMR grids, here we prefer to compute the exact values of $d$ on the entire grid, starting from the set of cells containing the points of $\pcloud$ and then proceeding by layers of neighbours, as described below in the so-called \textit{propagation} procedure (see \S~\ref{ssec:propagation}).

\subsection{Reinitialization}
\label{ssec:reinitialization}

A delicate issue when dealing with LSM is to guarantee that the level set function $\phi$ stays well-behaved during its evolution, hopefully close to a signed distance function, so that
\begin{equation}\label{eq:gradBound}
    \modulo{\nabla\phi} \approx 1.
\end{equation}
This is fundamental in order to get a watertight representation of the evolving surface and to exploit its geometrical properties in a reliable and stable way, since $\modulo{\nabla\phi}$ is well-defined and further away from zero. Also, for PDE models to be solved in $\Omega=\{ \phi\leq 0\}$, it is often crucial to be able to compute the normals to $\Gamma$ (i.e. by evaluating the gradient of $\phi$), at least close to the surface.

In this regard, even if we choose as initial data a signed distance function, the evolution governed by \eqref{eq:levelset:pde} will produce steep gradients in the solution and a so-called \textit{reinitialization} procedure needs to be performed. This constitutes a classical ingredient when dealing with LSM even if we point out that some alternatives to avoid the reinitialization step can be found in literature (see for example \cite{reini:1996,reini:2005,reini:2012,reini:FaPaTo:20}).

The most common way to reinitialize a level set function is by solving the steady-state HJ-type equation introduced in \cite{SuSmeOsh:1994} given by
\begin{equation}\label{eq:AMR:reinieq}
    \begin{cases}
        \phi_\tau + S(\tilde\phi)(\modulo{\nabla \phi}-1) =0,\\
        \phi(\vec{x},0) = \phi_0(\vec{x}) = \tilde\phi(\vec{x},t),
    \end{cases}
\end{equation}
where $\tilde\phi$ represents the possibly perturbed level set function to be reinitialized. In \cite{AMR:MiGuiBuGi:16} the authors solve the above equation following the explicit finite differences scheme described in \cite{AMR:MinGi:07}, to get a second-order accurate LSM on non-graded adaptive Cartesian grids.

Here instead we refer to the work of Saye \cite{AMR:Saye:14} in which the author proposes an efficient method for calculating high-order approximations of closest points on implicit surfaces that can be applied also on a general unstructured grid and in any number of spatial dimensions.

Actually, the procedure introduced in \cite{AMR:Saye:14} is applied just in the vicinity of the zero level set of $\phi$ and then, once we have reinitialized $\phi$ in these quadrants, we propagate the information all over $\grid$, as described in \S~\ref{ssec:propagation}.
\newline

Let us consider a fixed time $n$ of the evolution of the level set $\phi$ for which we have obtained the updated level set function $\tilde\phi$ that needs to be reinitialized. Along the lines of \cite{AMR:Saye:14}, we start by detecting the quadrants that contain the zero level set of $\tilde\phi$, namely the set $\grid_0\subset\grid$, as in \S~\ref{ssec:localization}. For each quadrant $j\in\grid_0$, a high-order closest point algorithm via Newton's method is applied as follows:
\begin{itemize}
    \item for each quadrant $j$ we create 4 points $\vec{x}_{j_s}$ in 2D (resp. 8 points in 3D) located at the centres of a subgrid $2\times2$ (resp $2\times2\times2$) of the quadrant $j$;
    \item each point $\vec{x}_{j_s}$ is projected onto the zero level set of the polynomial $\rec_j[\tilde\Phi]$, namely the reconstruction defined on the quadrant $j$ from the updated values $\{\tilde\phi_j\}_{j\in\grid}$;
    \item for each quadrant $j\in\grid_0$ we find its closest point $\vec{x}_j^*$ in $X=\{\vec{x}_{j_s}\}_{j\in\grid_0}$;
    \item the point $\vec{x}_j^*\in X$ is then used as the initial guess to compute the minimum distance of $\vec{x}_j$ from the zero level set of $\tilde\phi$ via a Newton's method that relies on the reconstruction associated to the point $\vec{x}_j^*$;
    \item finally, for each $\vec{x}_j$, the reinitialized value $\phi_j$ is set equal to this minimum distance, multiplied by the sign of $\tilde\phi_j$, namely preserving the sign of the level set function before the reinitialization.
\end{itemize}

In practice, the method described in \cite{AMR:Saye:14} is based on a two steps procedure to compute the closest point to the interface $\tilde\phi=0$ from a point $\vec{x}_j$ located in its vicinity, considering the zero level set of the polynomials $\rec_j$. In what follows, the points $\vec{x}_j^*\in X$ are referred to as \textit{seed points}.

\subsection{Propagation}
\label{ssec:propagation}

\begin{figure*}[tpb]
\centering
\begin{minipage}{0.32\linewidth}
    \centering
    \includegraphics[width=1.\linewidth]{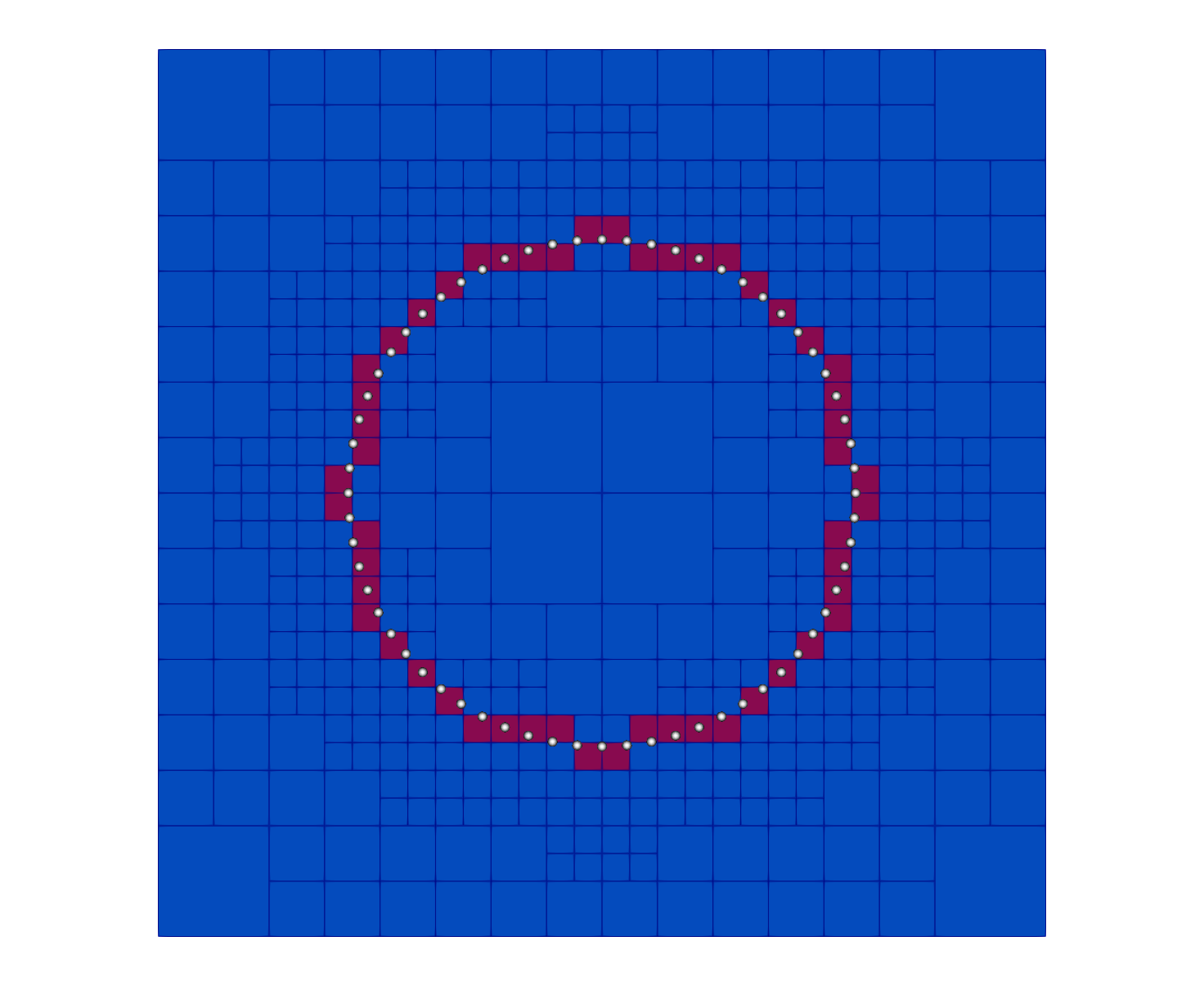}
\end{minipage}
\begin{minipage}{0.32\linewidth}
    \centering
    \includegraphics[width=1.\linewidth]{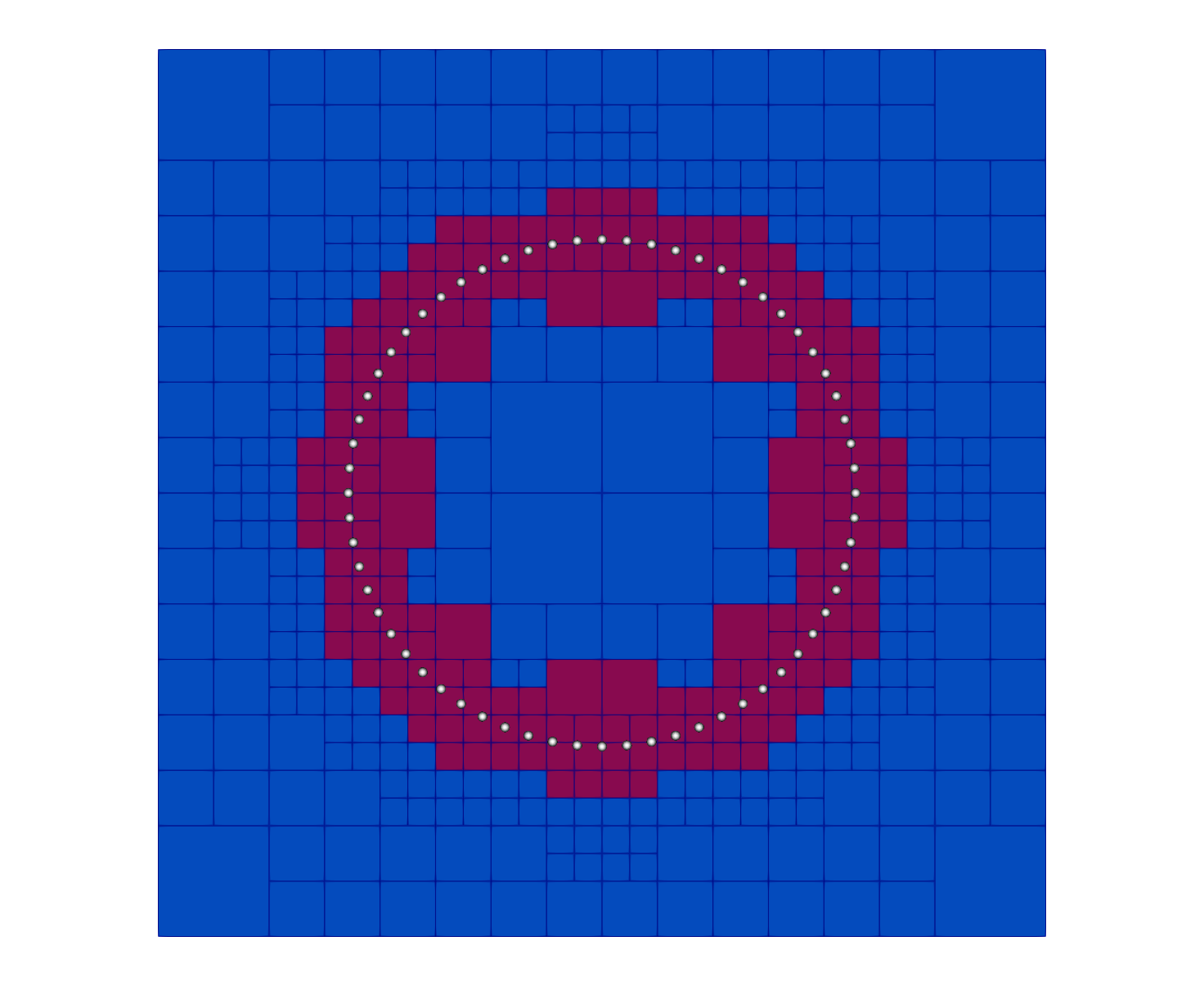}
\end{minipage}
\begin{minipage}{0.32\linewidth}
    \centering
    \includegraphics[width=1.\linewidth]{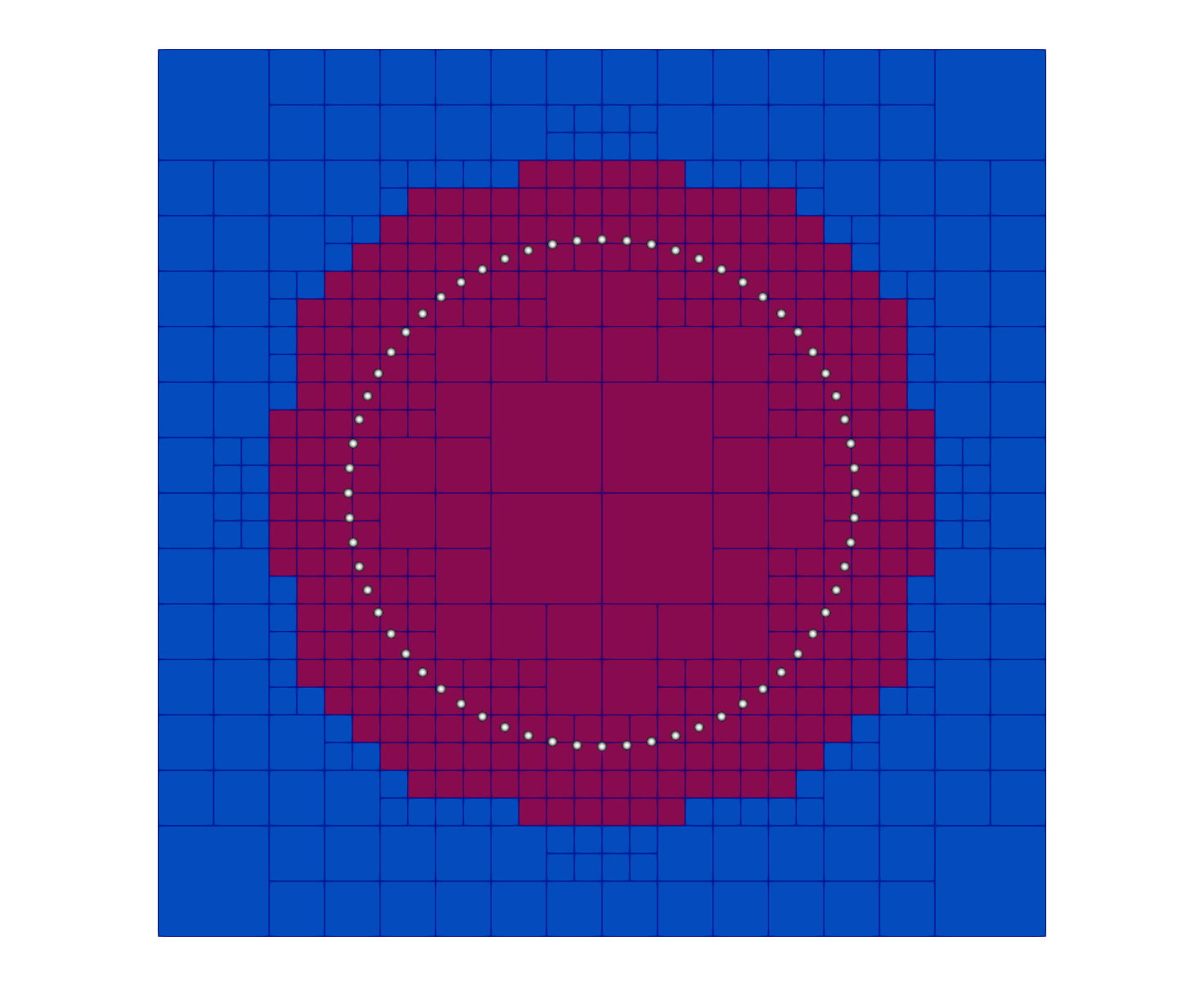}
\end{minipage}
\caption{Propagation of the distance function. Left: the initialized quadrants before the propagation step are the one containing one or some of the points in $\pcloud$. Center: quadrants with updated distance after one step of propagation. Right: quadrants with updated distance after two steps of propagation.} 
\label{fig:propagationAMR}	
\end{figure*}

Let us consider a subset $\doneSet^0\subset\grid$, on whose quadrants a certain quantity $g^0$ to be minimized is defined. In our specific case $g^0$ is the distance (resp. signed distance) of the quadrant centre $\vec{x}$ from a point $\vec{q}$, which is a member of $\pcloud$ (resp. a reinitialization seed point). We store this reference point as $\vec{q}_j$, attached to the quadrant $j$.

At the beginning of the propagation, data of quadrants in $\grid\setminus\doneSet^0$ are initialized to a high enough value. At the $m+1$-th iteration, let $\doneSet^m$ be the set of quadrants on which the value of $g$ has changed in the previous iteration, and $g^m$ be the current update of the function $g$ on the set $\cup_{k=0}^m\doneSet^k$. We will compute $g^{m+1}$ and $\doneSet^{m+1}$ as follows:
\begin{itemize}
    \item for each $j\in\doneSet^m$ we consider $\neighsFull^m_j$ the set of $j$ itself and its first neighbours, and we define the set $$\neighs^m\coloneqq\bigcup\limits_{j\in\doneSet^m} \neighsFull^m_j;$$
    note that in general $\neighs^m\cap\doneSet^m\neq\emptyset$, as we may need to further minimize the values updated in the previous propagation step, too;
    \item then, for each quadrant $i\in\neighs^m$ we compute the temporary value $\tilde g^{m+1}_i$ as
\begin{equation}
    \tilde g^{m+1}_i = \min_{k\in\doneSet^m}\modulo{\vec{x}_i - \vec{q}_k},
\end{equation}
where the point $\vec{q}_k$ is associated to $\vec{x}_k$;
\item we define the set of updated quadrants as
\begin{equation}\label{eq:doneSetIter}
    \doneSet^{m+1} = \{ i \in \neighs^m : \tilde g^{m+1}_i < g^m_i \}
\end{equation}
and set
\begin{equation}\label{eq:updatePropagIter}
g_i^{m+1}=
    \begin{cases}
       \tilde g^{m+1}_i \quad &\text{if} \quad i\in\doneSet^{m+1},\\
       g^m_i \quad &\text{elsewhere};\\
    \end{cases}
\end{equation}
\item if $\doneSet^{m+1}\neq\emptyset$ we iterate, otherwise we stop.
\end{itemize}

The procedure described above works exactly as it has been described, for the propagation of the distance function from the point cloud $\pcloud$. The quadrants containing the cloud points are initialized and then their information is propagated, including the already initialized quadrants in the propagation procedure. In Figure~\ref{fig:propagationAMR} three different phases of the propagation are shown ($m=0,1,2$). Red quadrants are the ones that have been initialized with the distance to a point $\vec{q}\in\pcloud$, if $m=0$, or updated lowering their value of the distance to $\pcloud$, for $m=1,2$; the blue ones are the ones in which the information has not arrived yet.

Regarding reinitialization, the procedure is done analogously, but in the definition of $\neighs^m$, $m\geq 0$, we exclude the quadrants in $\doneSet^0$ since we do not want their values to be changed by the propagation process. Also, the condition in \eqref{eq:doneSetIter} and \eqref{eq:updatePropagIter} are modified respectively as
\begin{equation}\label{eq:doneSetIterReini}
    \doneSet^{m+1} = \{ i \in \neighs^m :  \modulo{\tilde g^{m+1}_i} < \modulo{g^m_i} \}
\end{equation}
and
\begin{equation}\label{eq:updatePropagIterReini}
g_i^{m+1}=
    \begin{cases}
       S(g^0_i) \, \tilde g^{m+1}_i \quad &\text{if} \quad i\in\doneSet^{m+1},\\
       g^m_{i} \quad &\text{elsewhere},
    \end{cases}
\end{equation}
where $S$ is the sign function.

\paragraph{Remark} Of course, the propagation algorithm must be implemented within the parallel setup of the rest of the code. This requires to iterate the scheme above more than one time in order to allow the information to reach all the processors, namely the whole computational domain. Once the local set $\doneSet^m$ is empty, local boundary quadrants spread their information to the halo region in order to minimize the function $g$ in ghost quadrants. At this point, a communication occurs. Each processor collects possible minimizing values from its neighbours and checks if the function $g$ should be actually updated in its local boundary quadrants. If so, the set $\doneSet^m$ would be filled by the updated quadrants, being not any more empty, causing the local propagation iterations to start again. Otherwise, if there are no local updates, the corresponding processor will wait for the others to finish the current round of local propagation. When all the $\doneSet^m$ remain empty, the overall propagation procedure stops.

\subsection{Energy functional}
\label{ssec:energy}

To apply the SL scheme \eqref{eq:SL2d} with $p>1$ we finally need to approximate, at each time step, the value of the energy functional $E_p(\phi)$ defined in \eqref{eq:energy:ls}. 

We approximate the Dirac $\delta$ function by restricting the integration domain to the subset 
\begin{equation}\label{eq:zeroSet}
\grid_0 \coloneqq \{ j\in\grid \, : \, \modulo{\{i\in\neighs_j : \phi_i\phi_j\leq 0\} } > 0 \},    
\end{equation}
namely the set composed by the cells in $\grid$ where, at a specific time step, the front is located. Further, thanks to the reinitialization procedure,
we assume that the function $\phi$ is, at least locally around $\grid_0$, a signed distance so that $|\nabla\phi|=1$.

We thus approximate the energy as
\begin{equation}\label{eq:energyApprox}
    E_p(\phi) \approx \left( \sum_{j\in \grid_0} |d_j|^p (\dxmin)^{n-1} \right) ^{1/p},   
\end{equation}
where $n$ in \eqref{eq:energyApprox} states for the dimension of the problem and we assume that the refinement provides the maximum resolution for the quadrants in $\grid_0$.

\section{Numerical tests}
\label{sec:numericaltests}

In this section we present some representative numerical results in 2D and 3D, showing different features of the reconstruction process. Point clouds data will be both synthetic and non-synthetic ones coming from laser-scanning of real objects. 

The overall computation of the level set reconstruction is summarized in Algorithm~\ref{algo}. In a preliminary phase, one creates the computational grid, computes the distance function, sets the initial data, and adapts the grid to the first zero level set. Then, the evolution starts and the level set function is evolved until a final stopping criterion is satisfied. Differently from \cite{PreSe:25}, thanks to the fully adaptive framework, Algorithm~\ref{algo} is here performed just once, setting the desired maximum level of refinement $L$ from the beginning. Once convergence is achieved, some parameters and the reconstruction operator are changed in order to perform five last regularizing iterations to smoothen the final data. Moreover, in these final steps the adaptivity criterion is only based on the level set function $\phi$, and not on the distance function.

\begin{algorithm}
\caption{Adaptive level set reconstruction from a point cloud $\pcloud$}\label{algo}
\vspace{5pt}
\begin{enumerate}
\item Define the computational domain $\Omega'$ and create a quadtree/octree grid $\grid$ encompassing the dataset $\pcloud$.
\item \label{algo:d} Compute the distance function $d(\vec{x})$ at all grid points, see \S~\ref{ssec:distance}.
\item \label{algo:u0} Set initial data $\{\phi^0_j\}_{\vec{x}_j\in\grid}$ enclosing $\pcloud$.
\item \label{algo:firstAdapt} Cut and adapt using $\phi^0$ according to \S~\ref{ssec:localization}
 and \S~\ref{ssec:adaptivity}, respectively.
\item Loop:
  \begin{enumerate}
    \item \label{algo:banda} choose computational subgrid $\widetilde\grid\subset\grid$ as in \S~\ref{ssec:localization};
    \item \label{algo:energy} compute the energy functional \eqref{eq:energy:ls} as in \S~\ref{ssec:energy};
    \item \label{algo:sl} compute $\{\phi^{n+1}_j\}_{j\in\widetilde\grid}$ using the scheme defined in \S~\ref{sec:numericalscheme} coupled with the reconstruction operator from \S~\ref{ssec:P1_LS};
    \item \label{algo:reinit} reinitialize $\{\phi^{n+1}_j\}_{j\in\grid}$ as in \S~\ref{ssec:reinitialization} and cut according to \eqref{eq:cut};
    \item \label{algo:adapt} adapt using $\phi^{n+1}$ according to \S~\ref{ssec:adaptivity}.
  \end{enumerate}
\item The previous loop is repeated for five iterations with the reconstruction operator from \S~\ref{ssec:CW_LS} and different values of the parameters.
\end{enumerate}
\end{algorithm}

The datasets used for the tests were made available in the Digital Shape WorkBench of the AIM@SHAPE and VISIONAIR  projects \cite{frog}, in the 3D Scanning Repository of the Stanford University \cite{standford} or by the work of Hoppe et al. \cite{hoppe1992}. All the codes for the simulations have been implemented in C++ using the P4EST library \cite{AMR:P4EST:11} for grid management and MPI parallelization. The 3D tests have been performed on the cluster Galileo 100 hosted at CINECA\footnote{https://www.hpc.cineca.it/systems/hardware/galileo100/}, exploiting the resources assigned to ISCRA-C Projects\footnote{\textit{Adaptive Mesh Refinement in Level Set Methods} (HP10C7HWOL) \\ \textit{High-order schemes on Adaptive Meshes} (HP10CO8NC7)}.

\paragraph{Spatio-temporal discretization}
To define the computational domain $\Omega'$ we first rescaled the cloud $\pcloud$ in a box $[-1,1]^n$ and then set $\Omega'$ as an enlarged cube $[-M,M]^n$, $M>1$ that encompasses all the points in $\pcloud$. The constant $M$ and the maximum level of refinement $L$ are computed in order to have
\begin{equation}\label{eq:dxminCloudSize}
    \dxmin = \frac{\Delta \Omega'}{2^L} \approx C_\pcloud h_\pcloud,
\end{equation}
where $h_\pcloud$ denotes an estimate of the resolution of the cloud $\pcloud$ and $C_\pcloud$ is a constant set by the user. The value of $h_\pcloud$ is approximated by randomly choosing a sample made up of the $10\%$ of the points in $\pcloud$ and then computing the average of their distances from the cloud itself. Regarding $C_\pcloud$, having in mind the further applications to PDE models, we usually choose values between $0.25$ and $0.5$ in order to get a final reconstruction defined on a finer grid, compared to the cloud resolution.

Concerning the time step, for each quadrant $j\in\grid$, we set
\begin{equation}\label{eq:dt}
    \dt_j = \frac{3}{2}\dx_j,
\end{equation}
which clearly states the advantages of using the SL approach: to update the solution we are allowed to set $\dt = \Ogrande(\dx)$, not being prohibitively constrained by the parabolic term. In \eqref{eq:dt}, note that $\dx_j$ corresponds to the edge length of the actual active quadrant $j$.

\paragraph{Initial data} 
Regarding the initial data, differently from \cite{Zhao:2000}, this is simply set to be a signed distance function associated to a circle or a sphere containing all the data points. While this choice has the effect of slowing down the reconstruction process, it reduces the dependence on other parameters (see \cite{PreSe:25}), and prevents the risks arising from not evenly distributed datasets. In fact, especially in 3D, it is common to find datasets that present fake holes due to the supports used during the laser-scanning phase. These missing data could be misleading and need to be properly distinguished from actual topological features. That is why, exploiting the potential of adaptivity, in this work, we totally avoid these kind of issues choosing the simplest geometry and topology for the initial data. Moreover, in \S~\ref{ssec:tunnel}, we propose a modification to the level set equation \eqref{eq:levelset:pde} in order to capture diverse topologies. 

\paragraph{Stopping criterion} 
To stop the evolution, as we are looking for a minimum of the energy functional \eqref{eq:energy}, we resort to the following criteria: at time $n$, the algorithm stops if
\begin{equation}\label{eq:deltaEn}%
    \Delta_E^n = \frac{|\Bar{e}^k_{n-1}-\Bar{e}^k_{n}|}{\Bar{e}^k_{n}} < 10^{-4},
\end{equation}
where
\begin{equation}
    \Bar{e}^k_{n} = \frac{1}{k}\sum_{i=n-k+1}^n E_2(\phi^i)
\end{equation}
or after a maximum of $100$ iterations. Note that in \eqref{eq:deltaEn}, the energy functional \eqref{eq:energy} is computed with $p=2$ and we force the algorithm to do at least 10 iterations by setting $k=\min(n,10)$. Condition \eqref{eq:deltaEn} is the one proposed in \cite{He:2019} and is in practice a way to detect stationary points or flat areas of the energy functional. The evaluation of $\modulo{\phi(\vec{q})}$, for any $\vec{q}\in\pcloud$, or $\modulo{\phi_t}$ could constitute possible alternatives, as suggested in \cite{Zhao:2000,HaMi:16}.

\paragraph{Parameters setting} 
The parameters $p$ and $\mu$ appearing in \eqref{eq:levelset:pde} need to be chosen. According to \cite{PreSe:25}, we set them initially as
\begin{equation}\label{eq:paramPCurvIni}
    p = 1 \quad \text{and} \quad \mu \in [0.05, 0.2],
\end{equation}
and perform the loop 5 in Algorithm~\ref{algo} until convergence is achieved, namely when \eqref{eq:deltaEn} is satisfied. Then, we update them as
\begin{equation}\label{eq:paramPCurvFin}
    p = 2 \quad \text{and} \quad \mu = 1,
\end{equation}
and perform five more iterations of the loop 6 in Algorithm~\ref{algo}.
We recall that $p=1$ grants for a faster evolution and it thus useful given our choice for the initial data, while $p=2$ is more effective in reaching a steady state for the evolution. Moreover, setting $\mu=0$ simplifies the update formulas \eqref{eq:SL2d} by nullifying the parabolic term while also preventing from loosing too much details, especially on coarse grids, at a price of having a more rugged surface as final result. On the other hand, higher values are suggested to smoothen the solution and to handle possible changes of topology. Unless otherwise specified we initially choose $\mu=0.2$. While in \cite{HaMi:16,PreSe:25} smaller values of $\mu$ were strongly required to handle the coarser Cartesian grids, they are no longer needed now. In the final iterations, we set $\mu=1$ in order to reproduce the original model and to regularize the final reconstruction. An example of the effect of the parameter $\mu$ will be given in \S~\ref{ssec:cubosfere}.

\paragraph{Reconstruction operator} 
Following the same spirit of parameter selection, we choose different reconstruction operators for the loops $5$ and $6$ in Algorithm~\ref{algo}. In particular, we expect the $\Pone$ operator to be cheaper, less accurate, but closer to a polygonal reconstruction faithfull to the points in $\pcloud$. On the other side, a high-order reconstruction should provide more accurate and regular results, being also capable to reproduce complex features of the shape. According to these considerations and to our focus on PDE computations, referring to Algorithm~\ref{algo}, we use the $\Pone$ operator of \S~\ref{ssec:P1_LS} in loop 5, while we switch to the $\CWENO$ one of \S~\ref{ssec:CW_LS} for the regularizing iterations of loop 6, to get a smoother final result. More details about this choice will be give in the square test in \S~\ref{ssec:squareTest}.

\paragraph{Results evaluation}
Finally, in order to evaluate the quality of the reconstruction we compute, alongside \eqref{eq:deltaEn}, 
the $L^1$-norm of the error $Err_1^n$, when the exact signed distance function $\phi^*$ is given, and the average of the error on the points of the cloud  
\begin{equation}\label{errCloud}
    Err_{\pcloud}^n = \frac{\sum\limits_{\vec{q}\in \pcloud} | \rec[\Phi^n]\left( \vec{q}\right)|} {|\pcloud|},
\end{equation}
to evaluate how much the final reconstruction is attached to the data. 

\subsection{Square 2D test}
\label{ssec:squareTest}

\begin{figure*}
\centering
\begin{minipage}{0.45\linewidth}
    \centering
    \includegraphics[width=1.\linewidth]{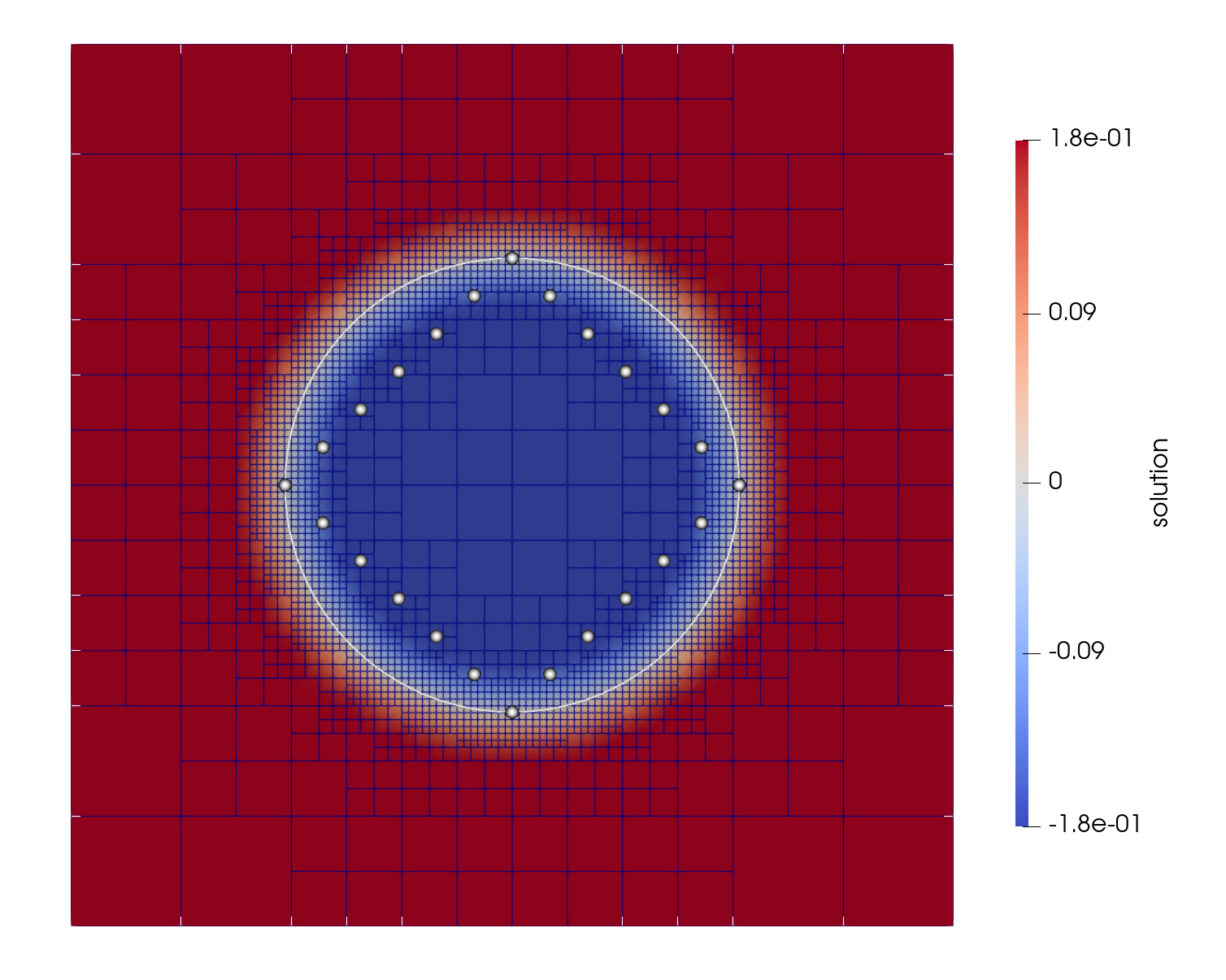}
\end{minipage}\hspace{0.2cm}
\begin{minipage}{0.45\linewidth}
    \centering    \includegraphics[width=1.\linewidth]{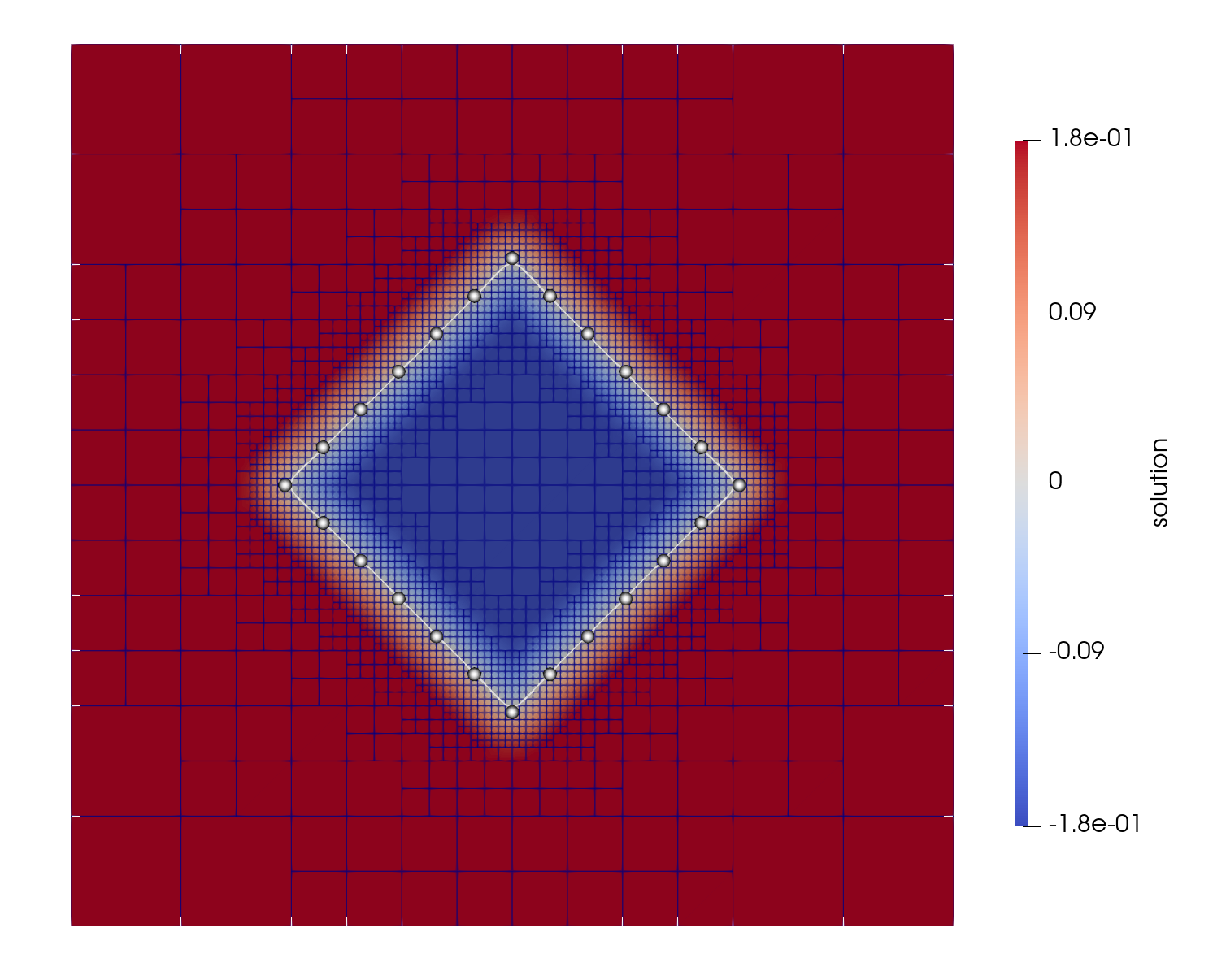}
\end{minipage}
    \caption{2D square test: from left to right, the initial data and the final reconstruction. The zero isocontour of the level set function $\phi$ is always represented with a white line.}
    \label{fig:squareReconstruction}
\end{figure*}

\begin{table*}
\begin{center}
\small
\begin{tabular}{|c|c|c|c|}
\hline & L1-error & Cloud error & CPU time (secs) \\ \hline \hline
$\Pone$ & $\num{9.57e-3}$ & $\num{4.21e-3}$ & $\num{2.62}$ \\ \hline 
$\CWENO$ & $\num{7.52e-3}$ & $\num{4.66e-3}$ & $\num{15.52}$ \\ \hline 
$\Pone + \CWENO$ & $\num{8.63e-3}$ & $\num{5.33e-3}$ & $\num{4.30}$ \\ \hline 
\end{tabular}
\end{center}
\caption{Comparison between the use of different reconstruction operators. The $\Pone$+$\CW$ strategy shows the best result in terms of errors versus computational time.}
\label{tab:squareComparison}
\end{table*}

\begin{figure*}
\begin{minipage}{0.33\linewidth}
    \centering
    \includegraphics[width=1.\linewidth]{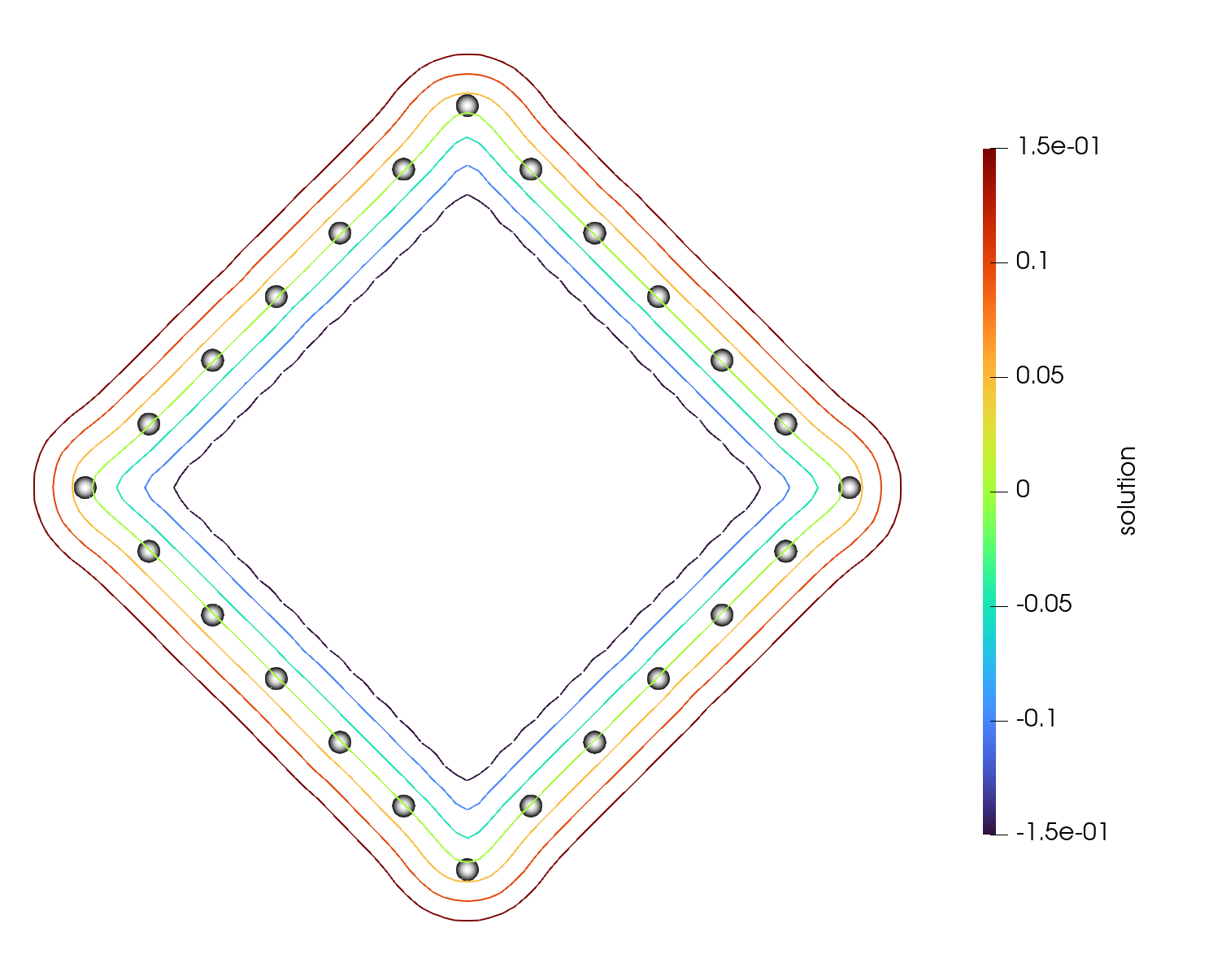}
\end{minipage}\hfill
\begin{minipage}{0.33\linewidth}
    \centering    \includegraphics[width=1.\linewidth]{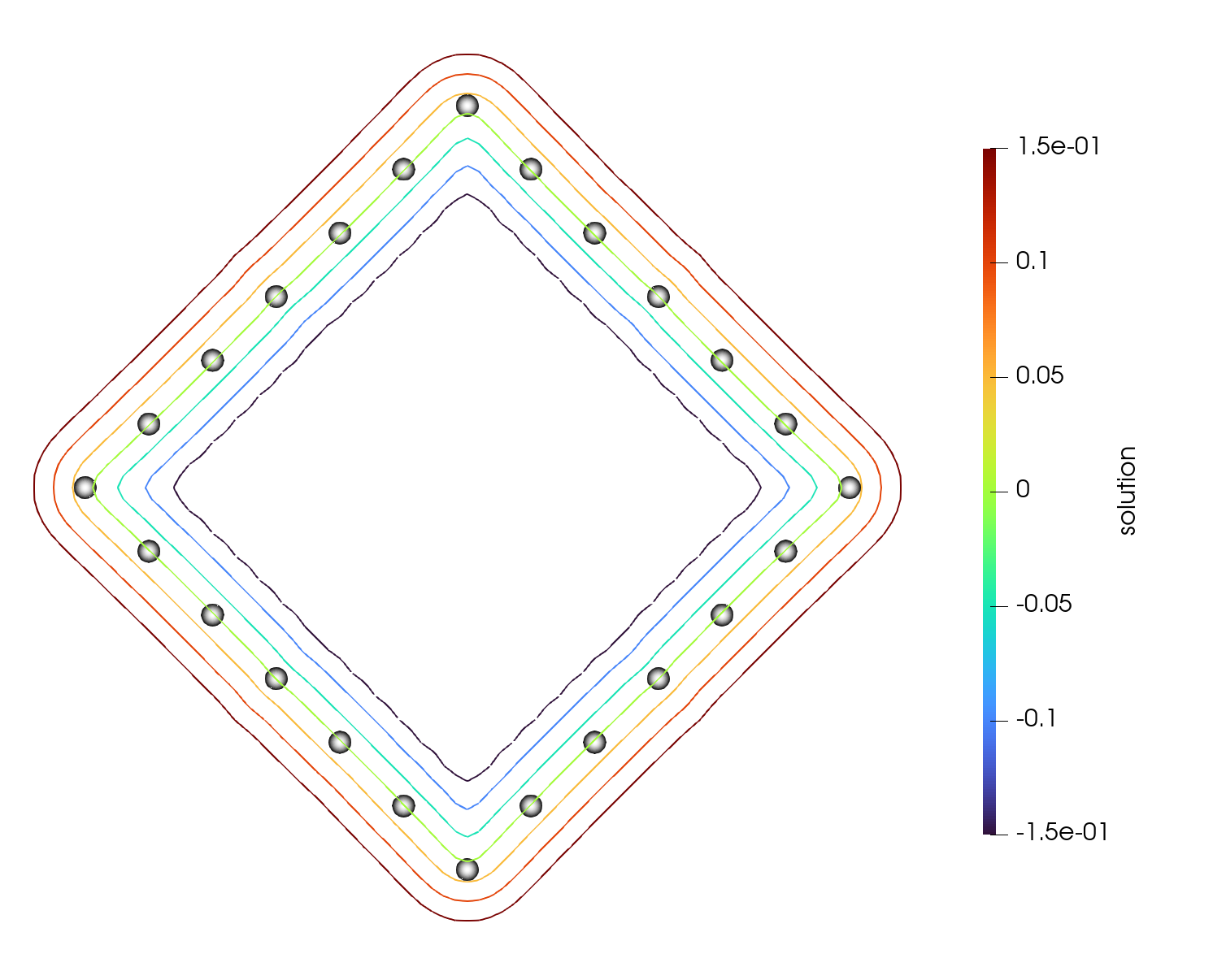}
\end{minipage}\hfill
\begin{minipage}{0.33\linewidth}
    \centering
    \includegraphics[width=1.\linewidth]{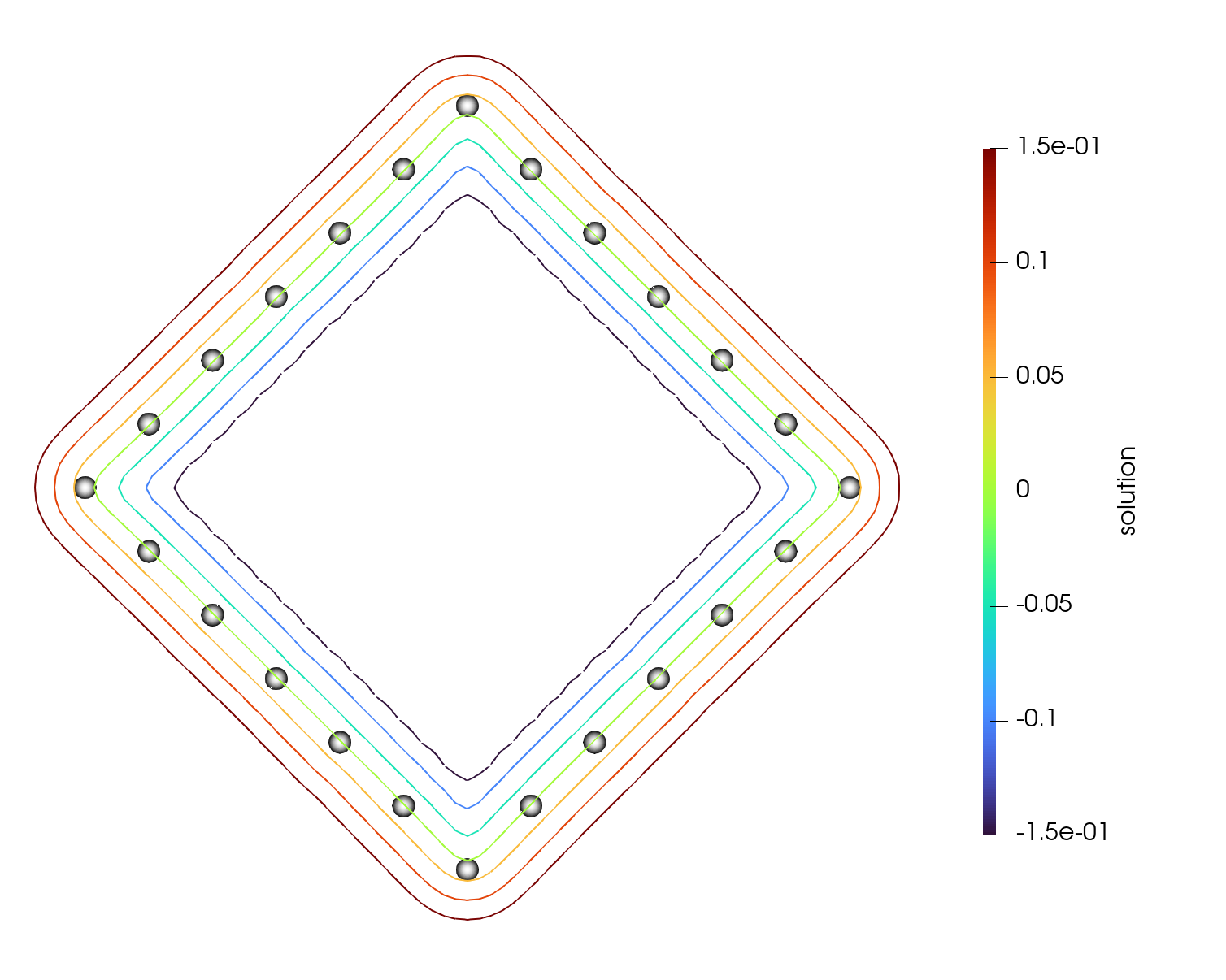}
\end{minipage}
    \caption{Final isocontours obtained by setting different reconstruction operators in the square test. From left to right: only $\Pone$, only $\CWENO$, the usual setting $\Pone + \CWENO$.}
    \label{fig:squareIsocontours}
\end{figure*}

We first consider the benchmark case of a square sampled by a point cloud $\pcloud$ made of 24 points whose resolution is $h_\pcloud\approx 0.23$ and for which the exact signed distance function $\phi^*$ is known. The initial data and the final reconstruction, obtained with the setting described above ($C_\pcloud=0.125$, $\dxmin=0.03$ and initial value $\mu=0.05$) are depicted in Figure~\ref{fig:squareReconstruction}, while in Table~\ref{tab:squareComparison} we show a comparison between strategies that employ different reconstruction operators: only $\Pone$, only $\CWENO$ and the usual setup $\Pone + \CWENO$, taking respectively $39$, $38$ and $39$ iterations to get the final result. One can notice that mixing the two operators, $\Pone$ and $\CWENO$, gives the best result in terms of accuracy versus computational time, even if we pay something in terms of fidelity to the dataset $\pcloud$. Similar considerations can be deduced from Figure~\ref{fig:squareIsocontours} where some isocontours of the final level set are depicted. We point out that using only the $\Pone$ operator the isocontours are less orthogonal around the vertices. On the other side, resorting to the $\CWENO$ reconstruction leads to a more accurate result, but Figure~\ref{fig:squareIsocontours} shows that there is no significant gain in applying it from the beginning.

\subsection{Heart 2D test}

\begin{figure*}
\centering
\begin{minipage}{0.45\linewidth}
    \centering
    \includegraphics[width=1.\linewidth]{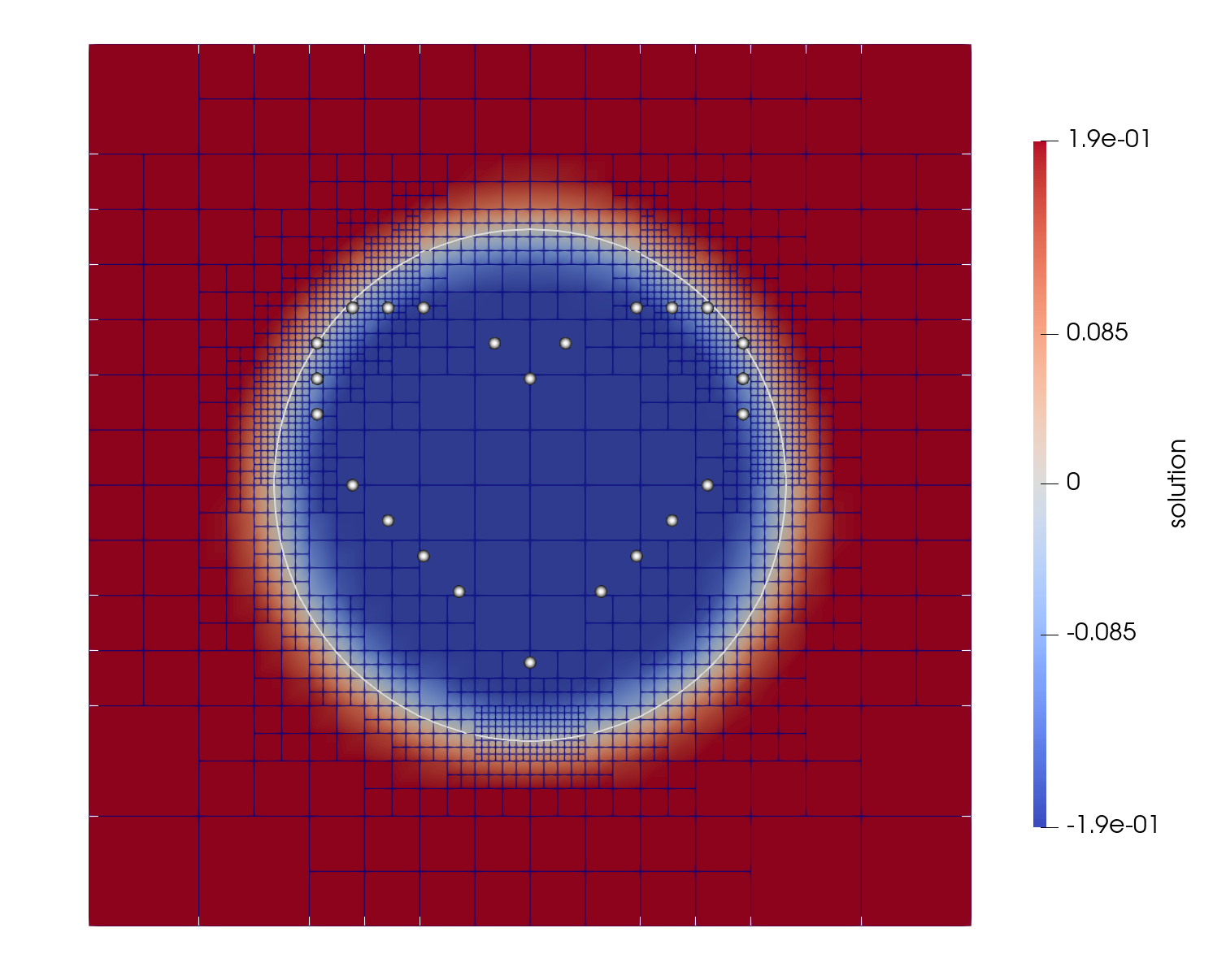}
\end{minipage}\hspace{0.2cm}
\begin{minipage}{0.45\linewidth}
    \centering    \includegraphics[width=1.\linewidth]{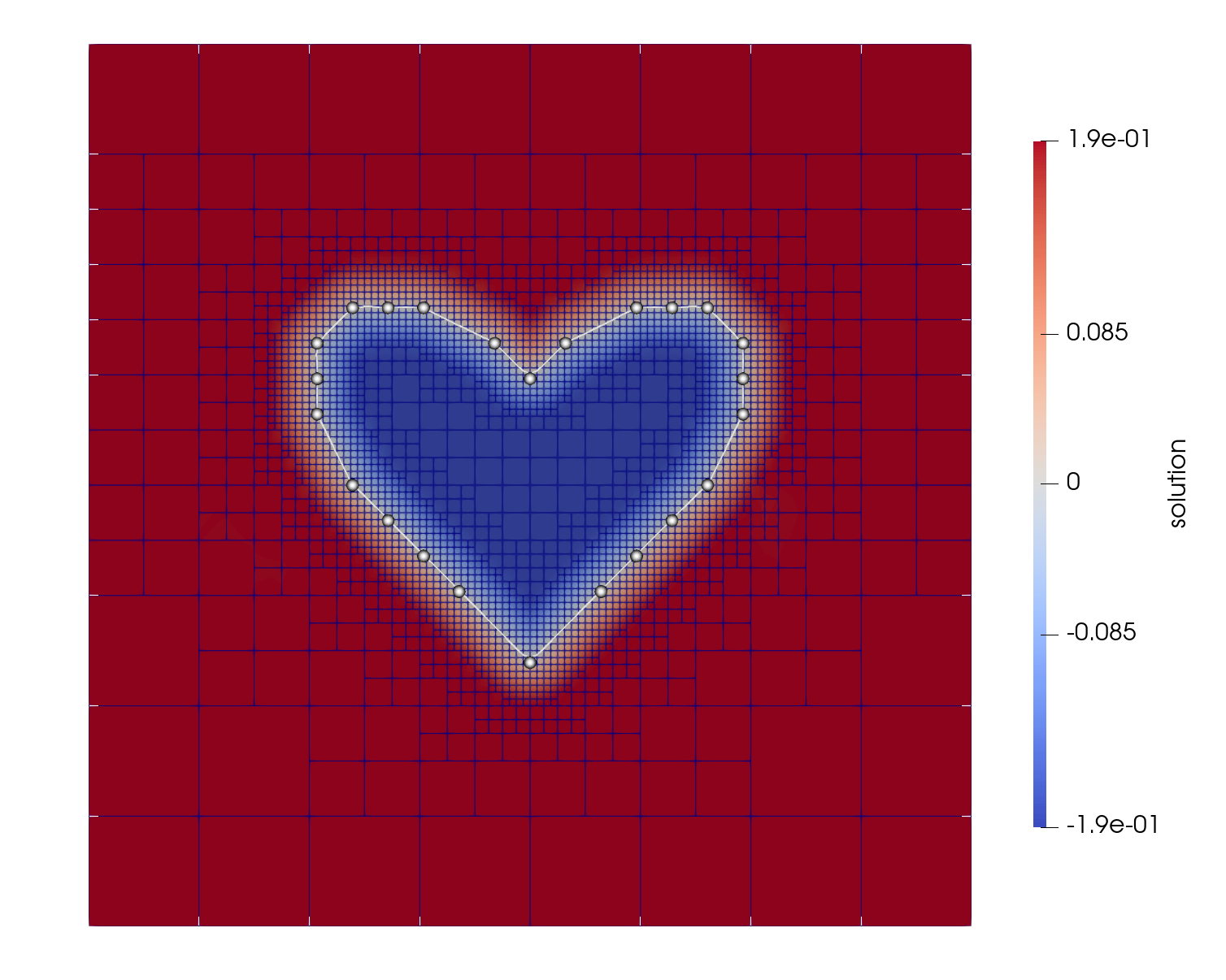}
\end{minipage}
    \caption{Reconstruction from a heart-shaped point cloud. From left to right: the initial data and the final signed distance function. The zero isocontour of $\phi$ is depicted by the white line.}
    \label{fig:heartReconstruction}
\end{figure*}

In the second test we consider a cloud of $24$ points representing a heart and having a resolution approximately equal to $\num{2.13e-1}$. Starting with $\dxmin=\num{3.24e-2}$, the reconstruction takes $43$ iterations and is capable to capture both the tips of the heart and its rounded parts. The initial data and the final reconstruction, with their zero isocontour, are shown in Figure~\ref{fig:heartReconstruction}.

\subsection{Teapot 2D test}

\begin{figure*}
\centering
\begin{minipage}{0.32\linewidth}
    \centering
    \includegraphics[width=1.1\linewidth]{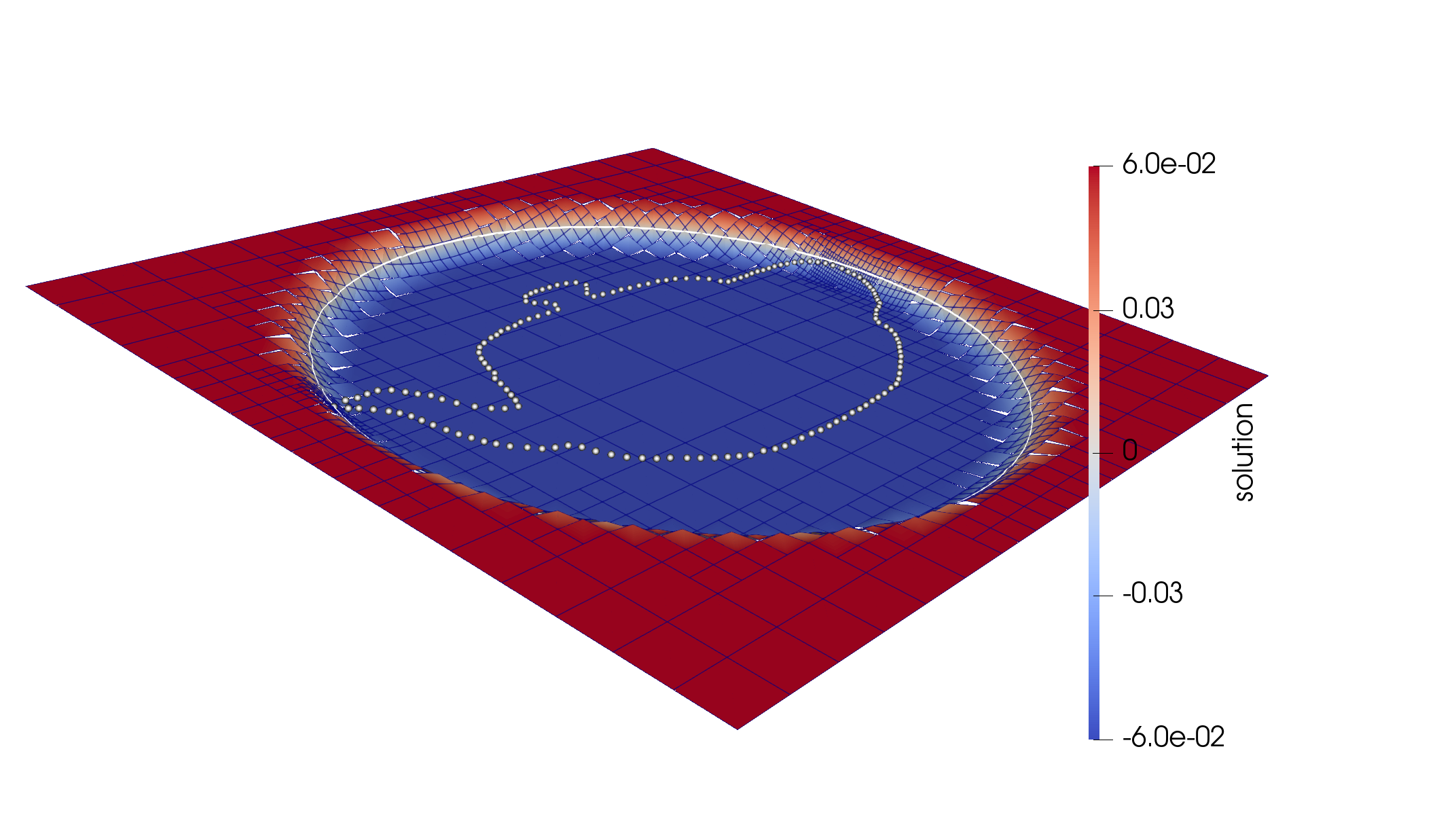}
\end{minipage}
\begin{minipage}{0.32\linewidth}
    \centering    \includegraphics[width=1.1\linewidth]{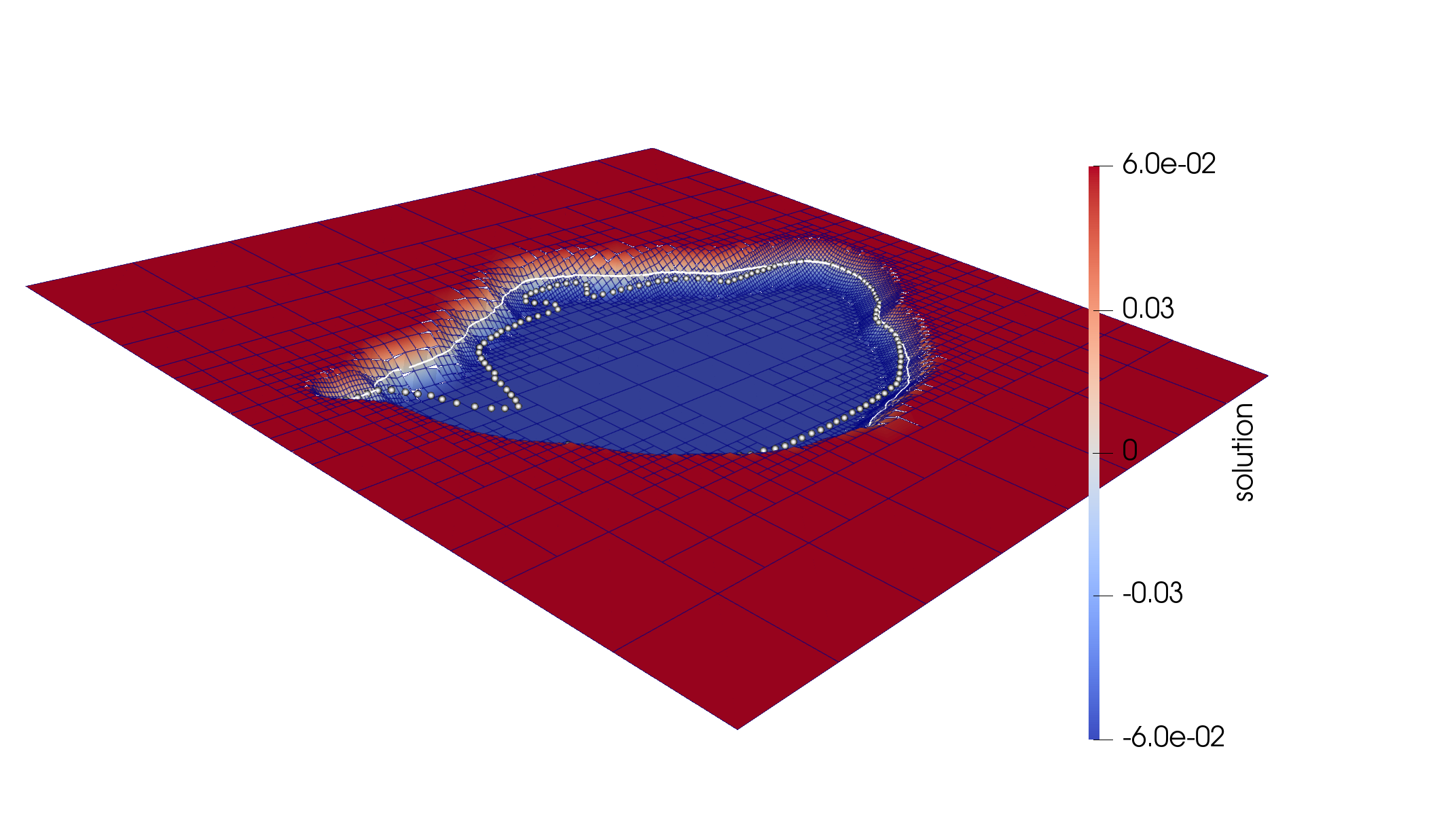}
\end{minipage}
\begin{minipage}{0.32\linewidth}
    \centering    \includegraphics[width=1.1\linewidth]{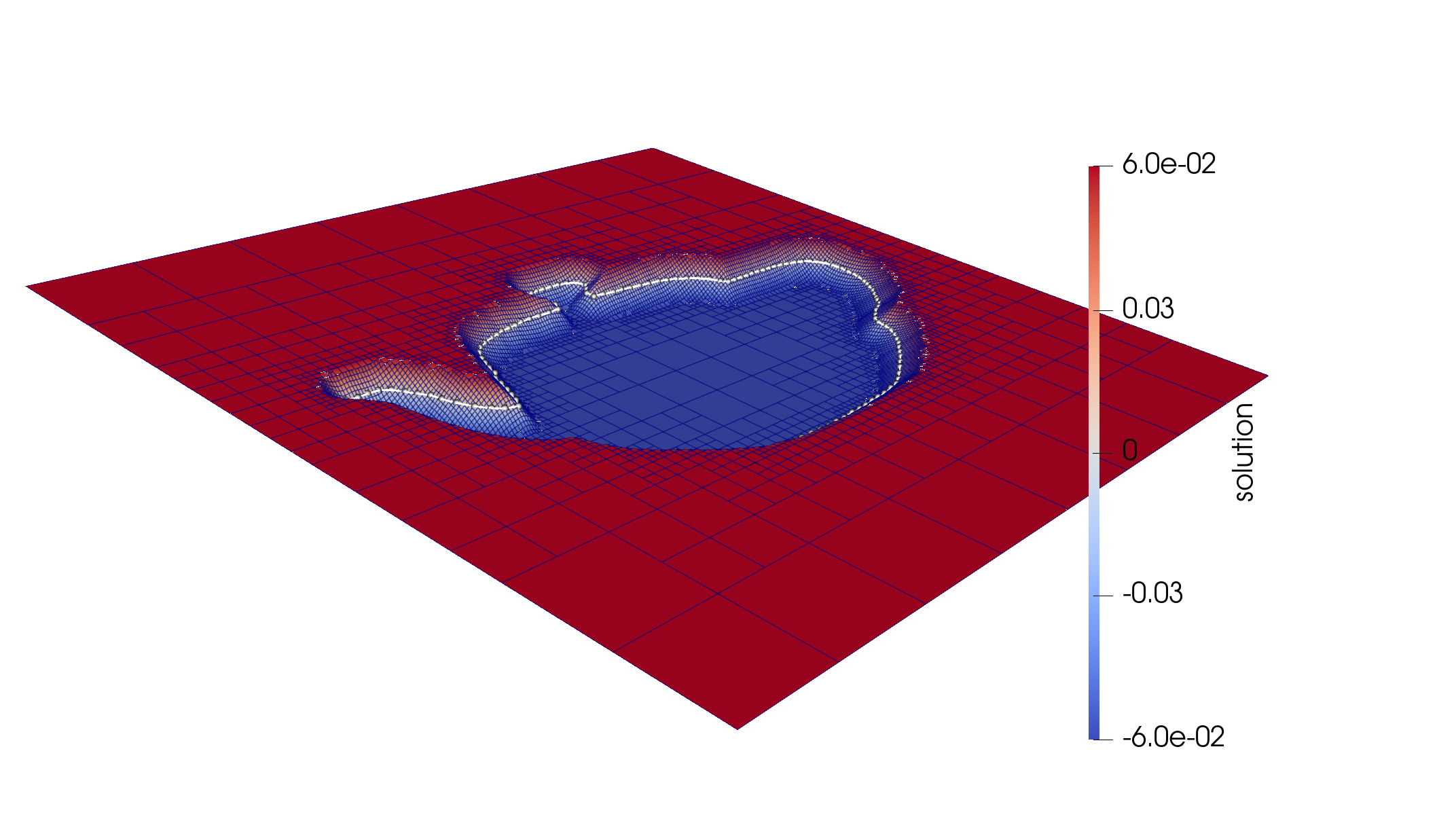}
\end{minipage}\hspace{0.2cm}
    \caption{Reconstruction from a teapot-shaped point cloud. From left to right: the initial data, an intermediate state and the final signed distance function. The zero isocontour of $\phi$ is depicted by the white line. Note that the grid around the white line is finer where the level set is closed to the point cloud.}
    \label{fig:tazzaReconstruction}
\end{figure*}

Here we consider a cloud of 189 points representing a teapot that presents different features, including rounded parts and fine details. In Figure~\ref{fig:tazzaReconstruction} we present three snapshots of the evolution of the level set and of the grid. One can appreciate that in the adaptive framework the resolution of the grid gradually increases as we get close to the data. The coarser grid allows to speed up the evolution when we are far from $\pcloud$, while the maximum level of refinement is effectively employed to capture all the details.

\subsection{``Cube\&Spheres'' synthetic 3D test}
\label{ssec:cubosfere}

\begin{figure}
\centering
\hspace*{-0.4cm}
\begin{minipage}{0.45\linewidth}
    \centering
    \includegraphics[width=1.3\linewidth]{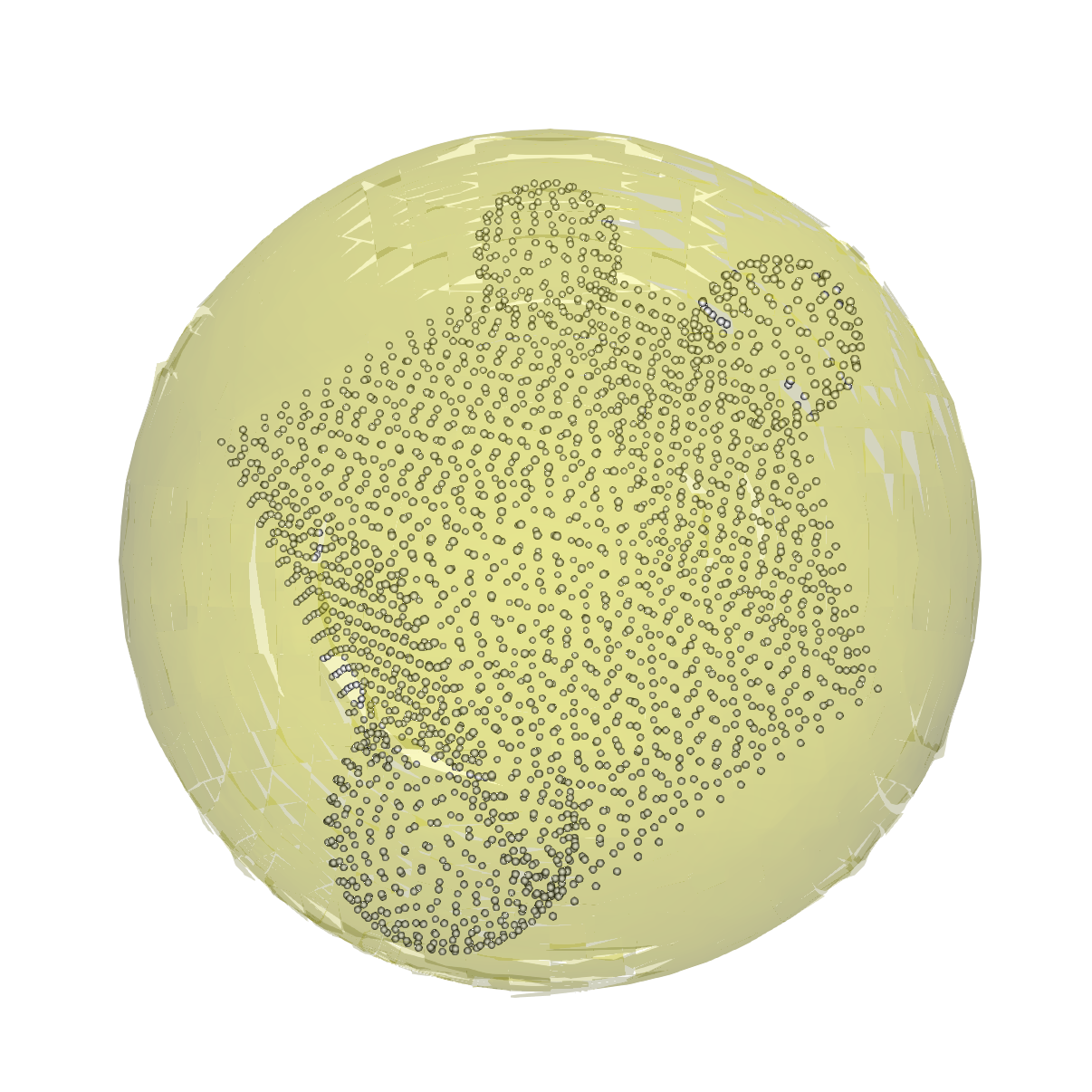}
\end{minipage}\hspace{0.1cm}
\begin{minipage}{0.45\linewidth}
    \centering    \includegraphics[width=1.3\linewidth]{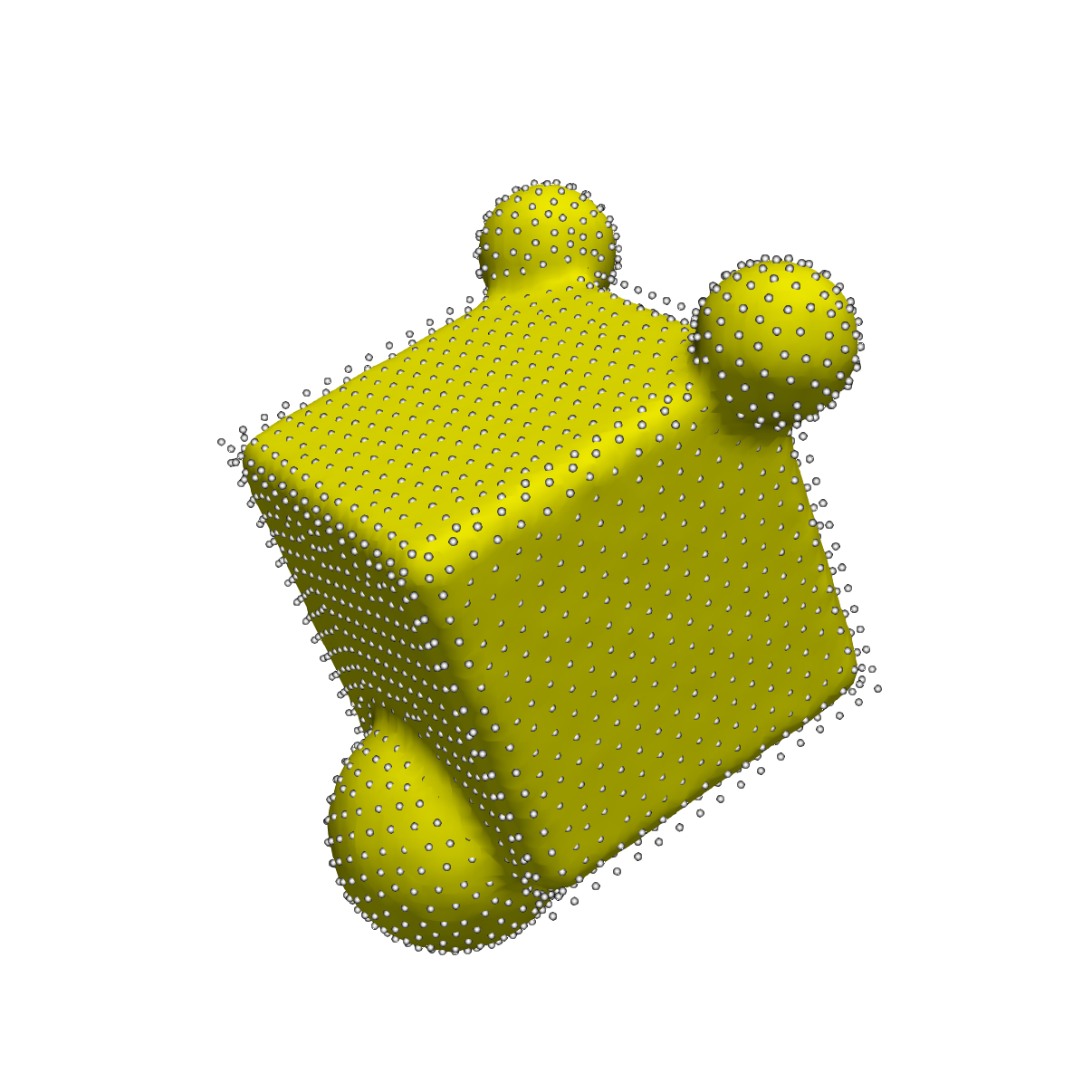}
\end{minipage}
    \caption{The initial surface and the reconstructed one of the ``Cube\&Spheres''.}
    \label{fig:cubosfereReconstruction}
\end{figure}

\begin{figure*}
\centering
\hspace*{0.5cm}
\begin{minipage}{0.45\linewidth}
    \centering
    \includegraphics[width=1.1\linewidth]{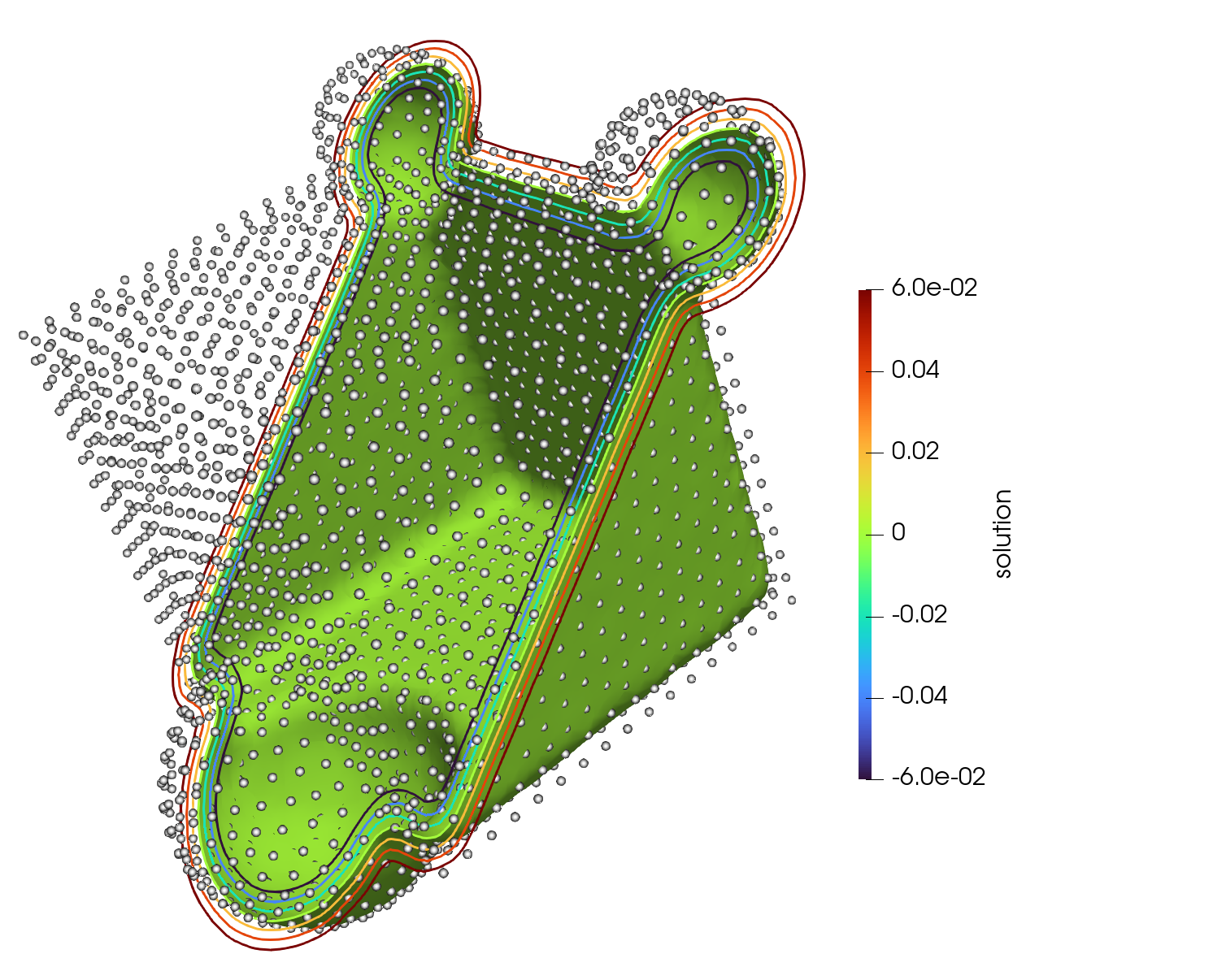}
\end{minipage}
\begin{minipage}{0.45\linewidth}
    \centering
    \includegraphics[width=1.1\linewidth]{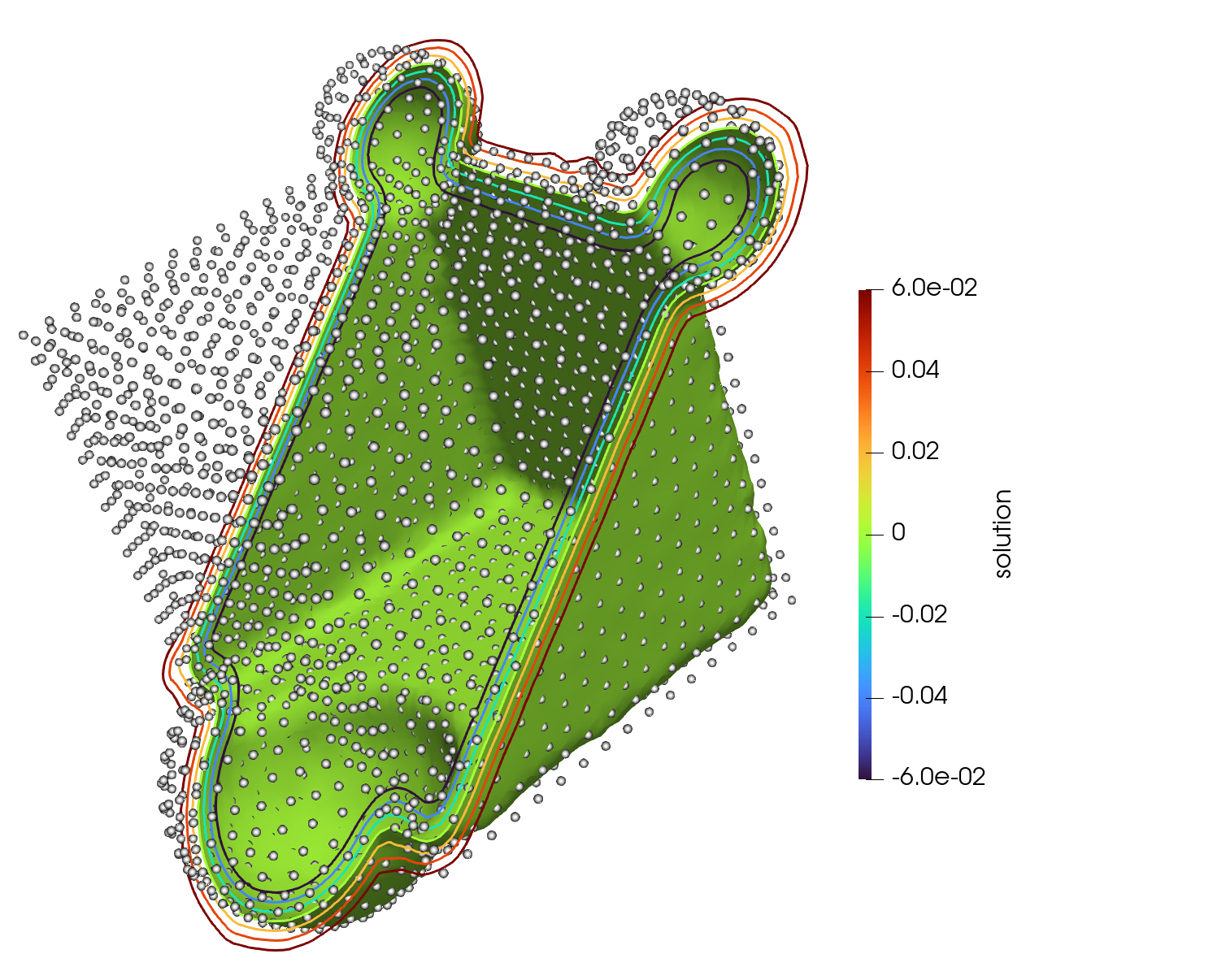}
\end{minipage}
    \caption{Evenly spaced isocontours computed from the ``Cube\&Spheres'' final level set for different initial values of $\mu$. On the left, $\mu=0.2$ gives a smoother result, whereas on the right, $\mu=0.05$ gives a more rugged surface.}
    \label{fig:cubosfereIsocontours}
\end{figure*}

We start the 3D session of tests with a domain composed by a cube joined with three spheres, sampled by 2346 points, which we refer to as ``Cube\&Spheres''. The cube edge length is $0.8$, the first sphere has radius $0.25$ and centre at the middle of an edge of the cube, while the other two had radius $0.15$ and were centred onto the two vertices of the opposite edge of the cube. The geometrical object was rotated in such a way that no face nor edge were aligned with the background grid. We find $h_\pcloud\approx \num{5.60e-2}$ and start the evolution with $\dxmin=\num{2.38e-2}$, by setting $C_\pcloud=0.5$. After $28$ iterations we get the reconstruction shown in Figure~\ref{fig:cubosfereReconstruction} and, since the exact solution is given, we can evaluate our approximation both in terms of $L^1$-error and cloud fidelity equal to $\num{4.32e-3}$ and $\num{4.45e-3}$, respectively. As we expected, higher values of $\mu$ lead to more regular solutions. This is clear from Figure~\ref{fig:cubosfereIsocontours}, where we show some isocontours of the final level set computed by setting initial values of $\mu$ equal to $0.2$ and $0.05$, respectively. At the expense of exact vanishing on the point cloud ($Err_{\pcloud}=\num{3.57e-3}$ for $\mu=0.05$), one can appreciate the improved regularity of the final reconstruction.

\subsection{Real data}
\label{ssec:realData}

\begin{figure*}
\vspace{-0.7cm}
\begin{minipage}{0.33\linewidth}
    \centering
    \includegraphics[width=1.2\linewidth]{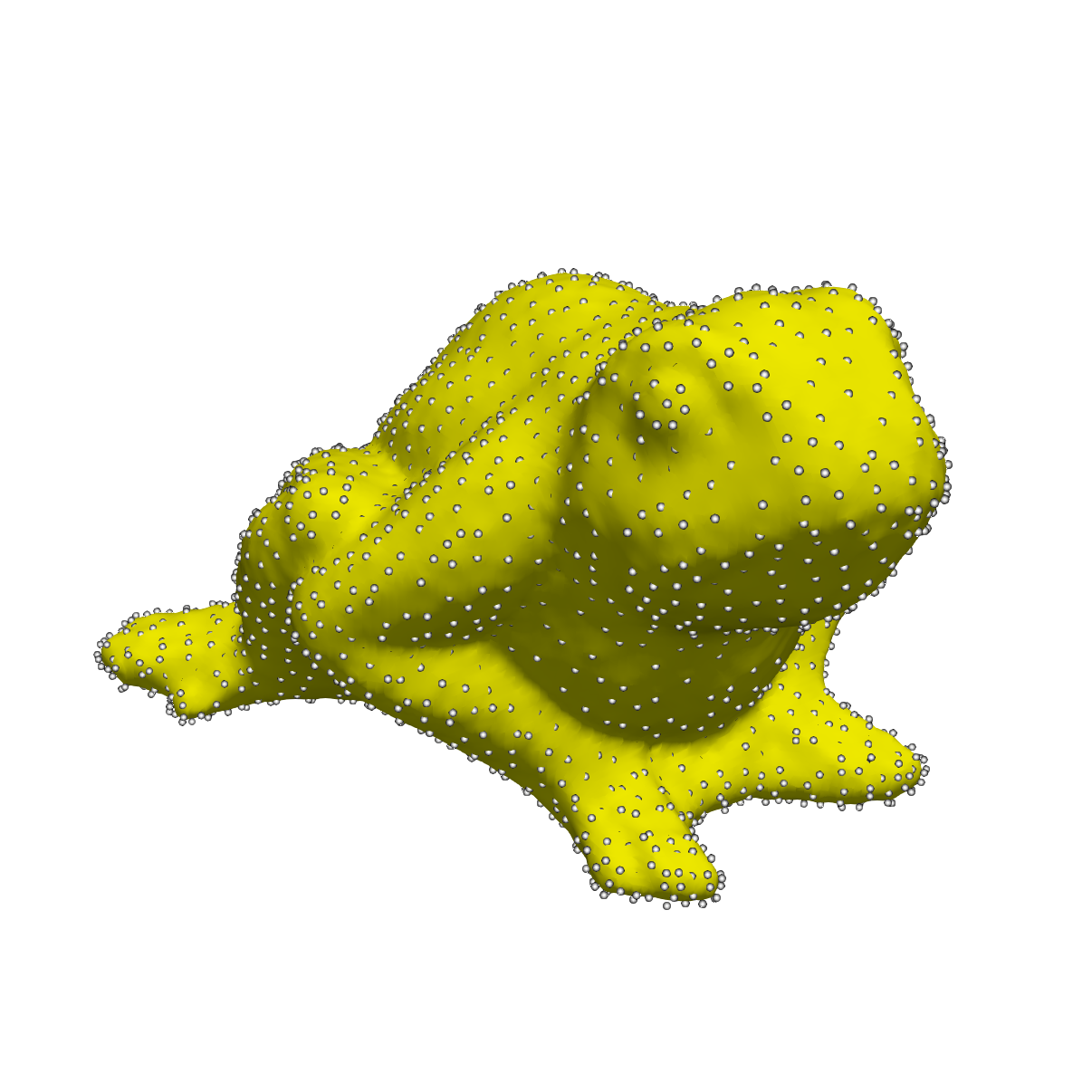}
\end{minipage}\hfill
\begin{minipage}{0.33\linewidth}
    \hspace{0.5cm}    
    \includegraphics[width=1.\linewidth]{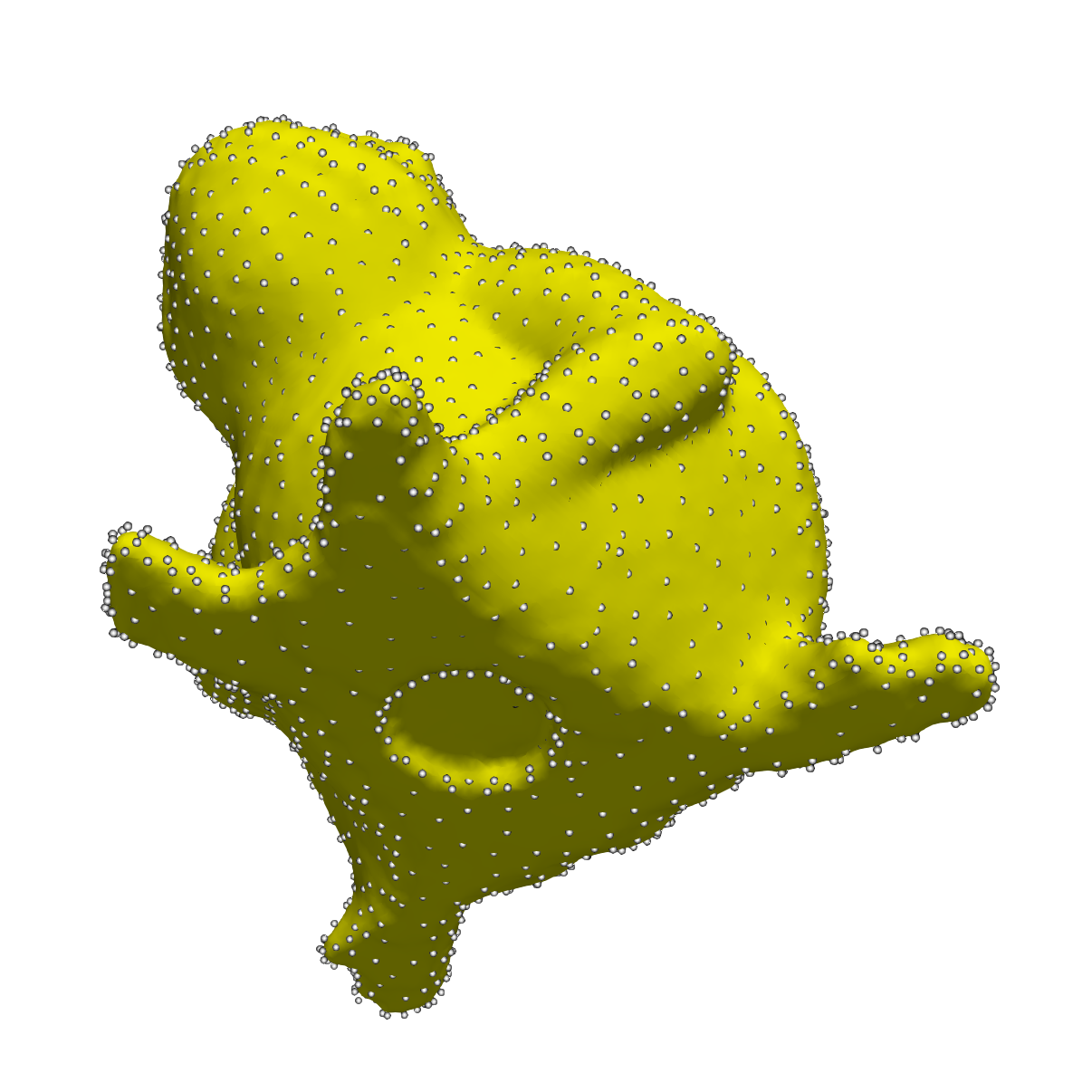}
\end{minipage}\hfill
\begin{minipage}{0.33\linewidth}
\vspace{0.3cm}
    \centering
    \includegraphics[width=0.73\linewidth]{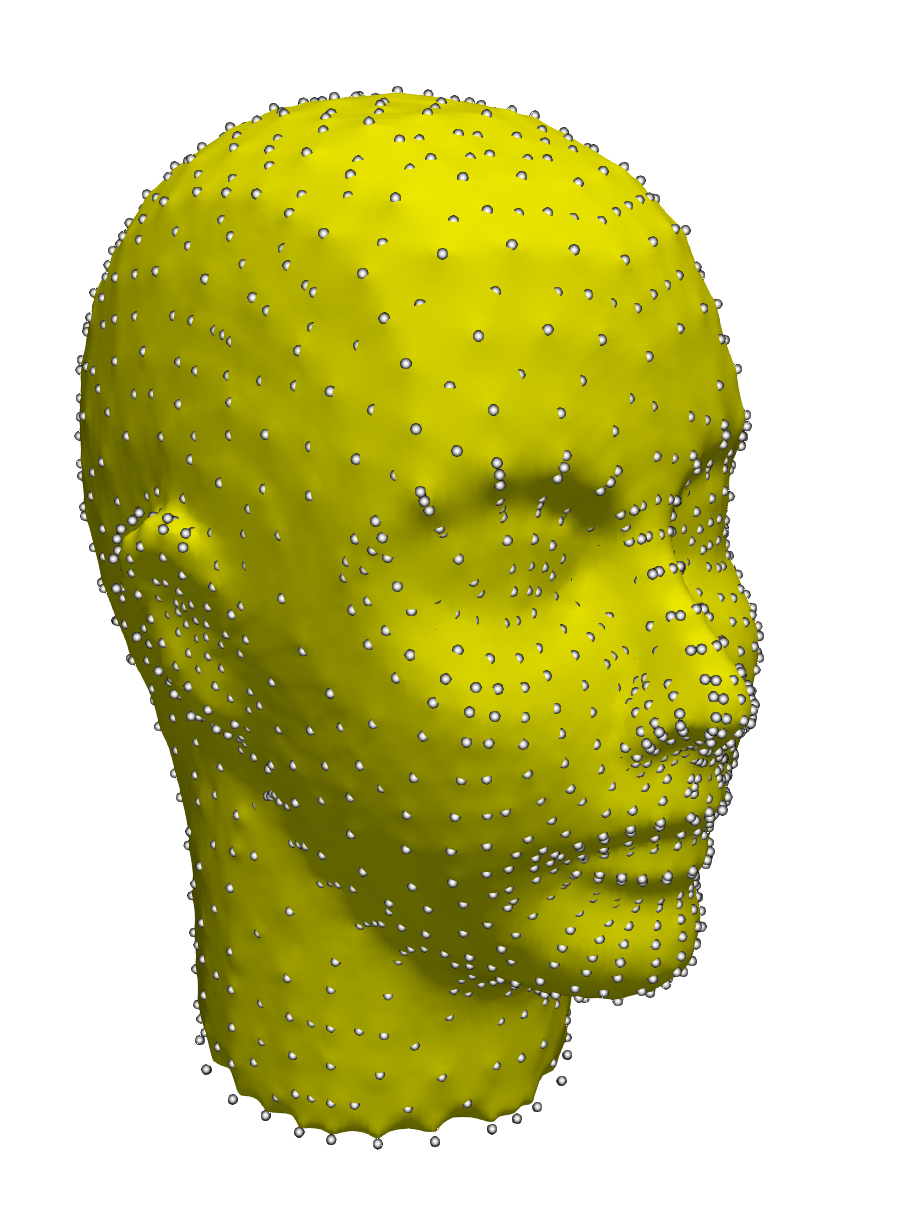}
\end{minipage}
    \caption{Reconstructed surfaces from the ``Frog'' and the ``Head'' datasets. Note how the level set is able to handle missing data at the bottom of the ``Frog'' cloud.}
    \label{fig:frogheadReconstruction}
\end{figure*}

Here we test our method using data sets coming from laser-scanning of real objects, namely the ``Frog'' and the ``Head'' point clouds coming respectively from \cite{frog} and \cite{hoppe1992}. Their reconstructions are shown in Figure~\ref{fig:frogheadReconstruction}. 

In the ``Frog'' test, we elaborate a cloud made up of $2512$ points with $h_\pcloud\approx\num{4.11e-2}$ and use $\dxmin=\num{1.43e-2}$, by setting $C_\pcloud=0.25$. The reconstruction takes $44$ iterations with a final $Err_\pcloud=\num{2.37e-3}$. This error is mostly affected by some parts of the reconstruction where the cross section of the shape is comparable to the resolution of the cloud, namely where the curvature is high or the thickness is very thin. Among the important advantages, we point out that this type of cloud has missing data that we would like to be plugged by our level set function, as it actually does (see Figure~\ref{fig:frogheadReconstruction}). 

On the other side, the ``Head'' cloud consists of $1881$ points having an estimated resolution $h_\pcloud\approx\num{3.83e-2}$. Setting $C_\pcloud=0.25$ we start with $\dxmin=\num{1.09e-2}$ and compute the reconstruction in $37$ total iterations getting a final error on the cloud $Err_\pcloud=\num{1.80e-3}$. 

\subsection{Complex shapes}
\label{ssec:complexShapes}

\begin{figure}
\hspace{-0.2cm}
\begin{minipage}{0.48\linewidth}
    \includegraphics[width=1.2\linewidth]{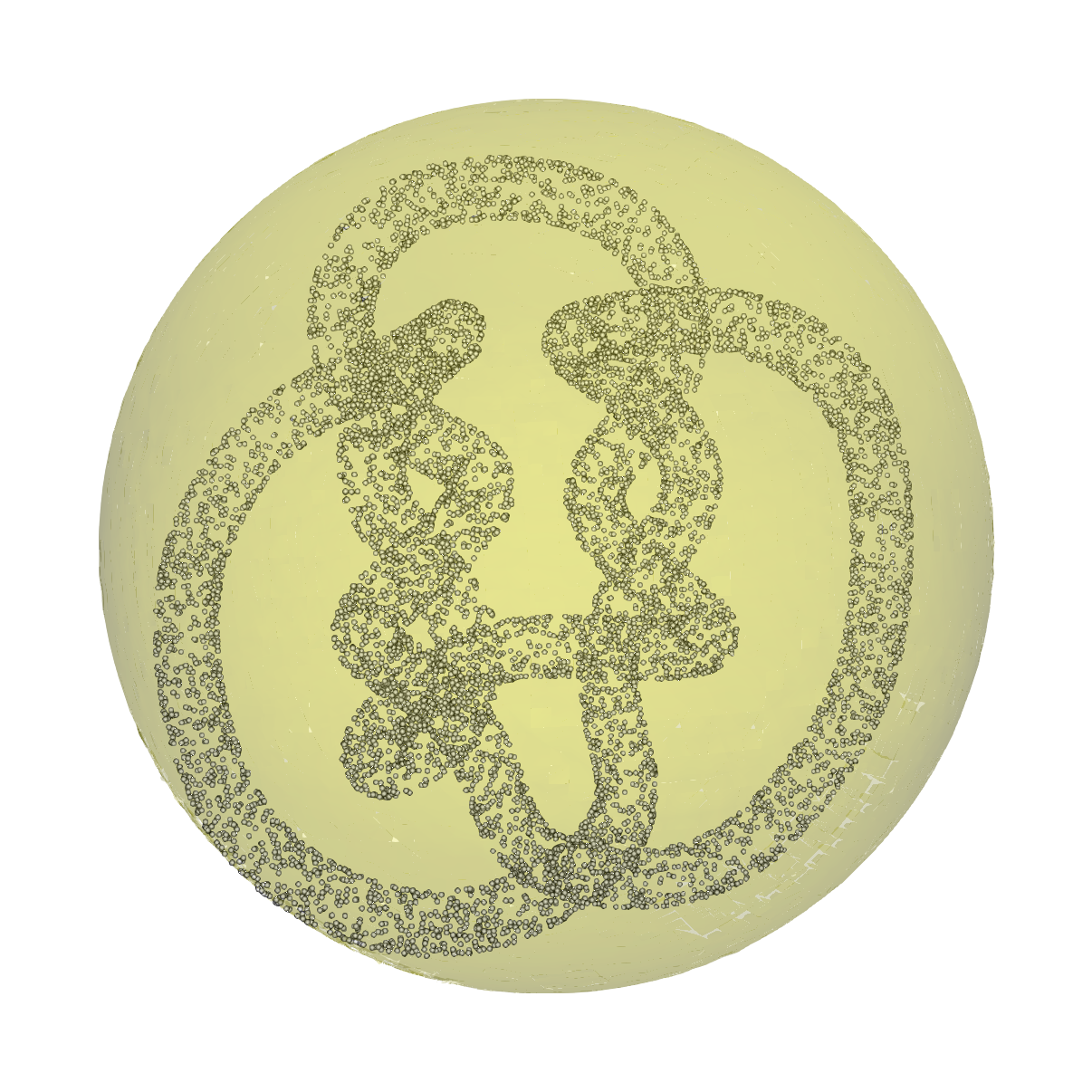}
\end{minipage}\hspace{0.1cm}
\begin{minipage}{0.48\linewidth}
    \includegraphics[width=1.2\linewidth]{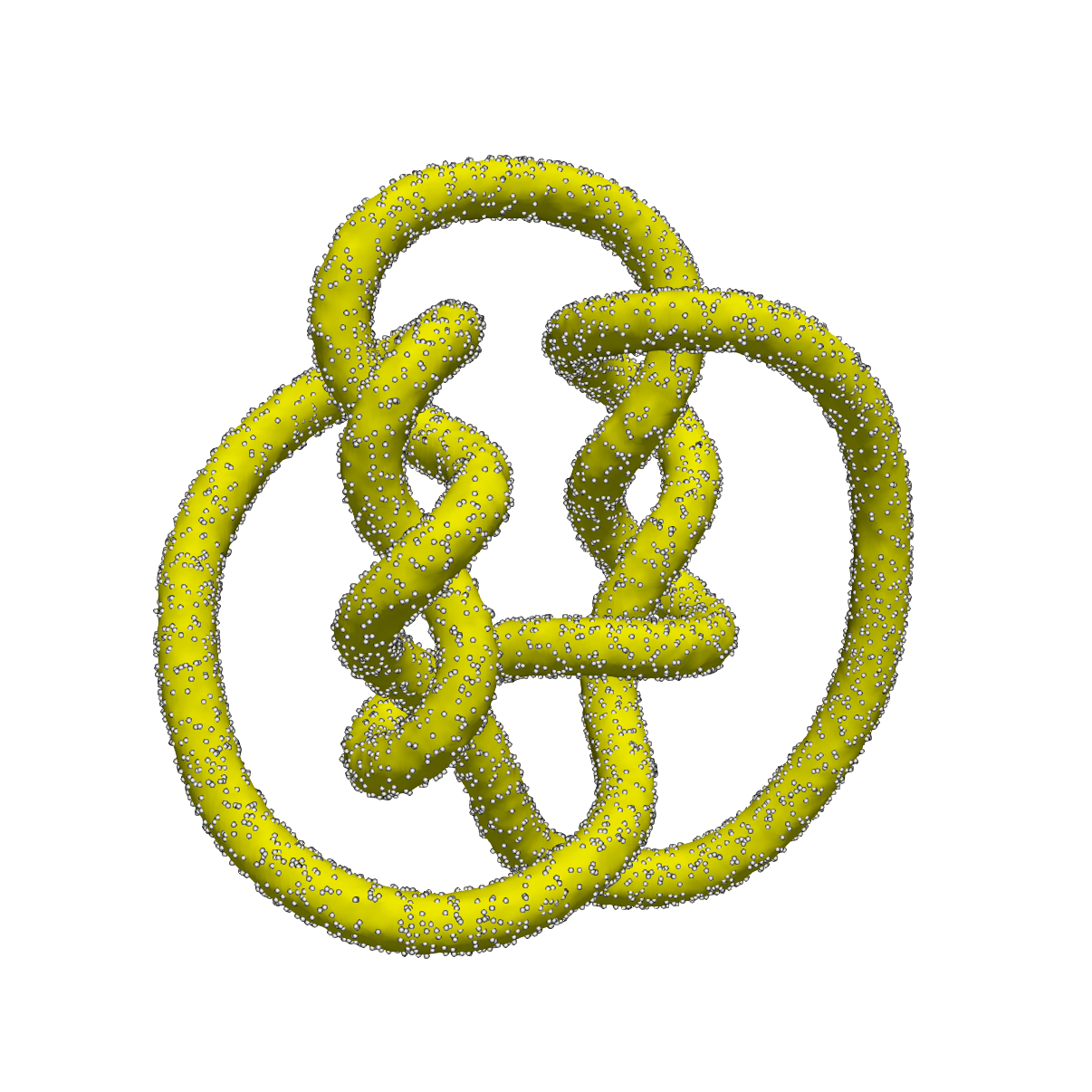}
\end{minipage}
    \caption{From left to right: the initial data and the final reconstruction for the ``Knot'' dataset. Note how the evolving level set is able to handle changes not only in geometry, but especially in topology.}
    \label{fig:knotReconstruction}
\end{figure}

In what follows we want to test the ability of our method to handle possible changes in topology. We thus consider a complex shape representing a knot, sampled by a point cloud of $10000$ points (\cite{hoppe1992}), having a resolution $h_\pcloud\approx\num{1.39e-2}$. From Figure~\ref{fig:knotReconstruction} one can notice that, even if we start the evolution with the simplest spherical initial data, we are able to capture all the features of the shape.

\subsection{Comparison with the Cartesian framework}

\begin{figure*}
\centering
\begin{minipage}{0.33\linewidth}
    \centering
    \includegraphics[width=1.1\linewidth]{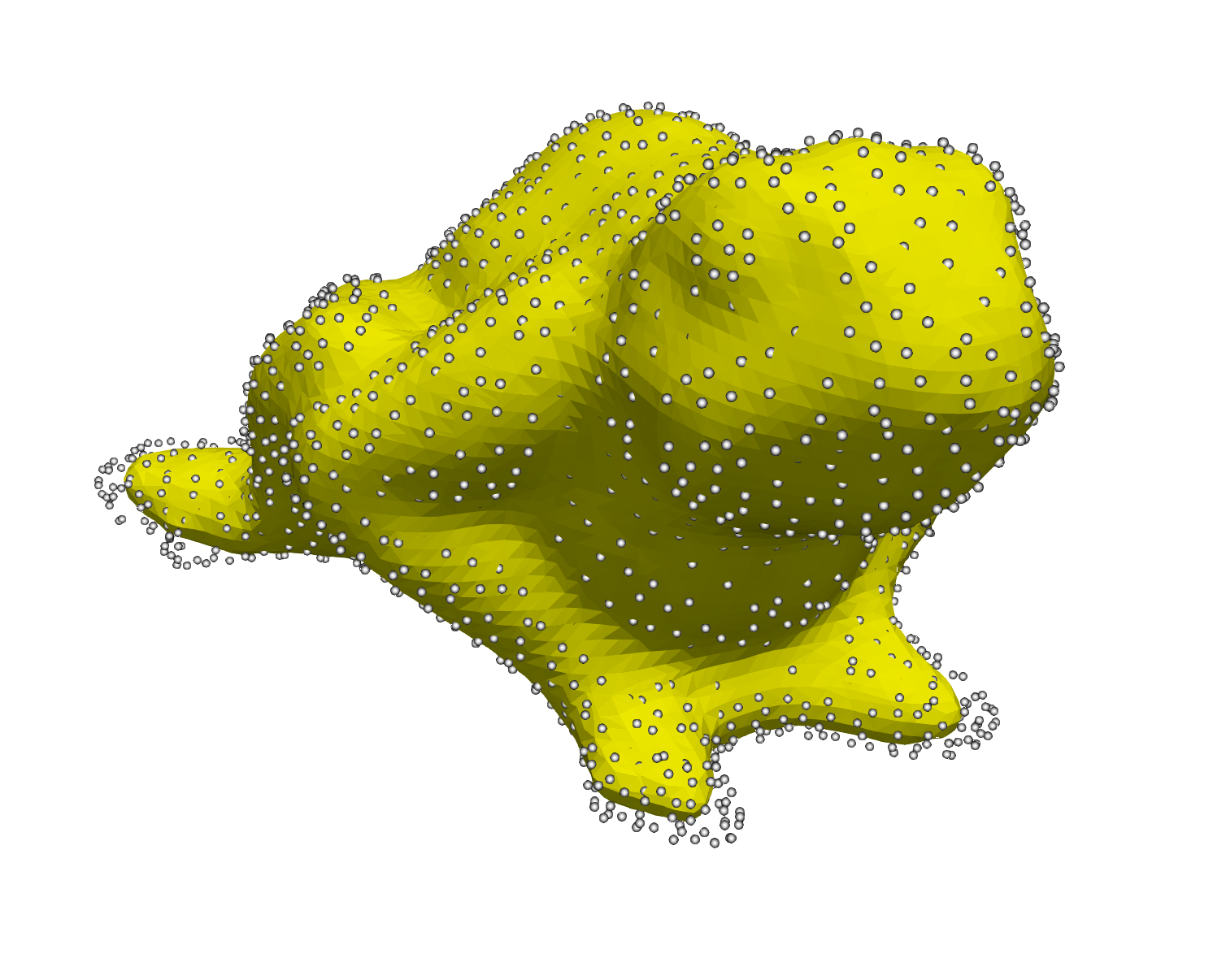}
\end{minipage}\hfill
\begin{minipage}{0.33\linewidth}
    \centering    \includegraphics[width=1.1\linewidth]{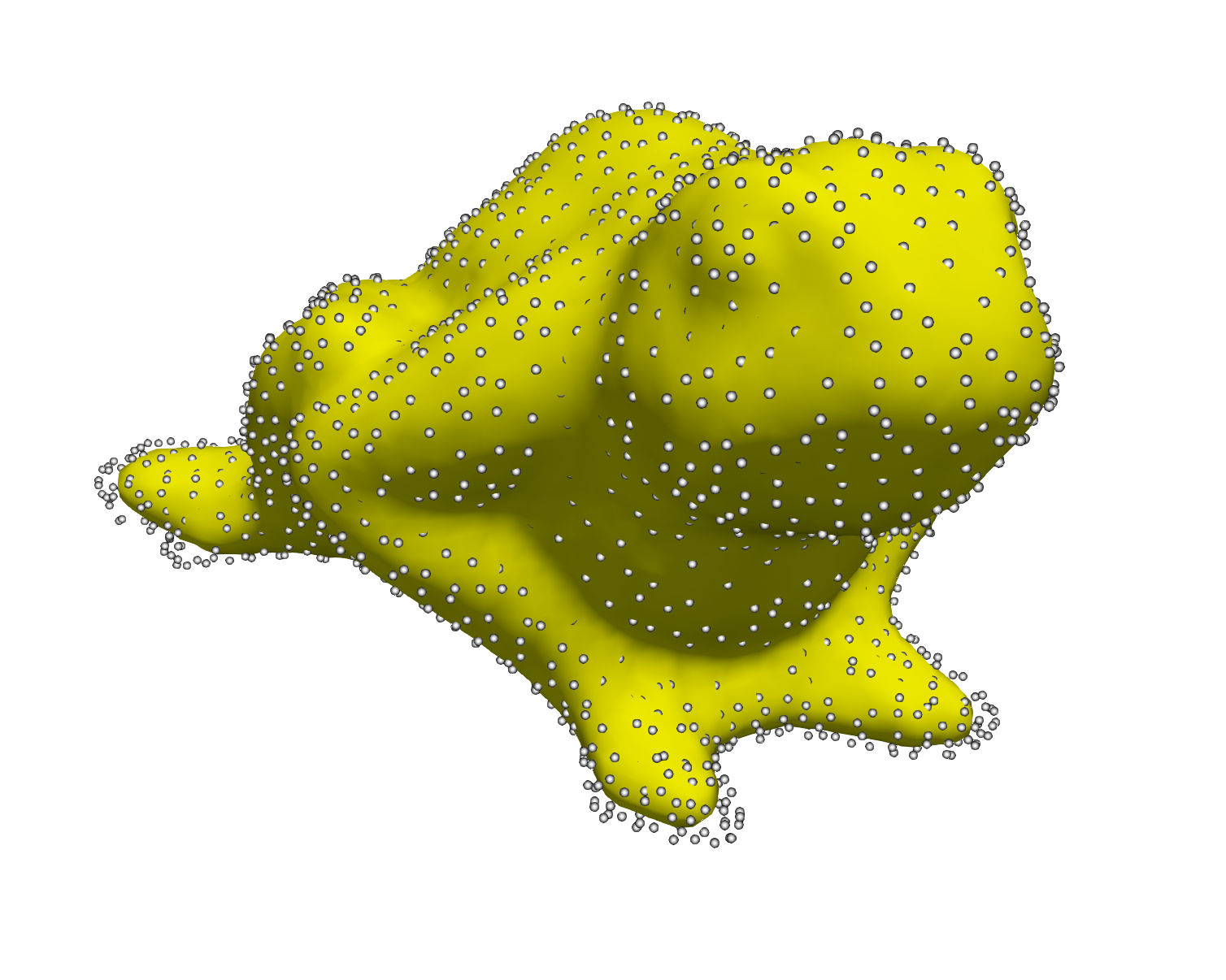}
\end{minipage}\hfill
\begin{minipage}{0.33\linewidth}
    \centering    \includegraphics[width=1.1\linewidth]{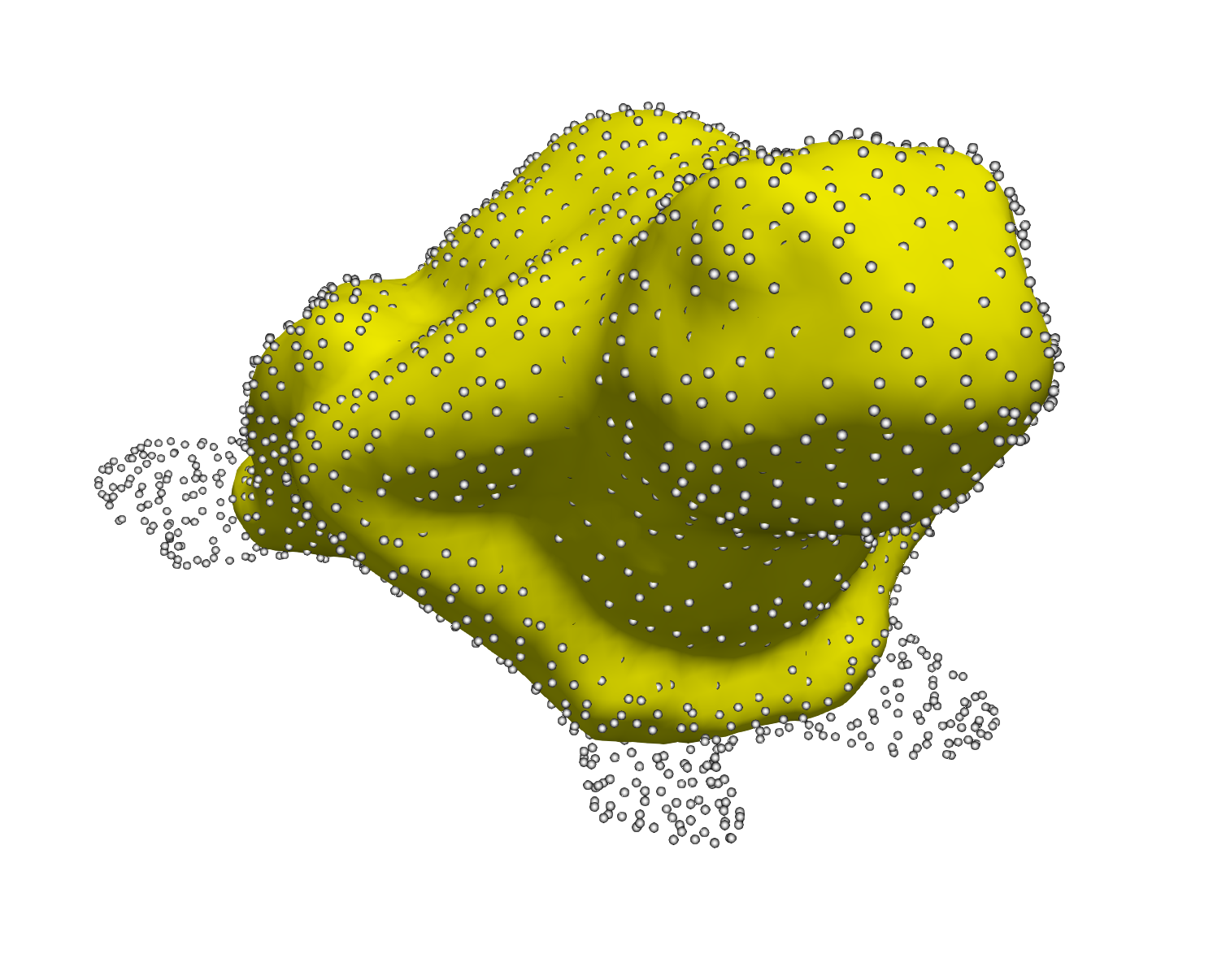}
\end{minipage}
    \caption{Comparison between the novel adaptive framework versus the Cartesian one of \cite{PreSe:25}, for the ``Frog'' dataset. From left to right: the final zero level set obtained with the adaptive method, the Cartesian one employing a single fine grid, and the Cartesian one employing two grids. While the first two results are comparable, in the last case, the evolution on the coarser grid is faster but completely fails to resolve the paws.}
    \label{fig:cartesianComparison}
\end{figure*}

\begin{table*}
\begin{center}
\small
\begin{tabular}{|c|c|c|c|c|}
\hline & Cloud error & Iterations & \makecell{CPU time\\(secs)} & \makecell{Degrees of\\freedom} \\ \hline \hline
Adaptive & $\num{7.65e-3}$ & $\num{28}$ & $\num{138}$ & $\num{59235}$ \\ \hline 
Cartesian (one grid) & $\num{6.75e-3}$ & $\num{29}$ & $\num{154}$ & $\num{599326}$ \\ \hline 
Cartesian (two grids) & $\num{2.16e-2}$ & $\num{66}$ & $\num{372}$ & $\num{409536}$ \\ \hline 
\end{tabular}
\end{center}
\caption{Comparison between the novel adaptive framework versus the Cartesian one of \cite{PreSe:25}, for the ``Frog'' dataset. The error on the cloud, the number of iterations, the total CPU time, and the number of degrees of freedom are considered.}
\label{tab:cartesianComparison}
\end{table*}

In this section we want to investigate the improvements achieved by this novel adaptive approach, compared to the previous Cartesian one, presented in \cite{PreSe:25}. The comparison is made by considering the ``Frog'' dataset of \S~\ref{ssec:realData} and performing the simulation on a laptop equipped with an Intel Core Ultra 9 185H CPU and 32 GB of RAM, employing $14$ ranks for the computations.

The adaptive reconstruction has been obtained with $\dxmin=\num{3.42e-2}$, with the usual setting described in this work and it is depicted in the left panel of Figure~\ref{fig:cartesianComparison}. When applying the Cartesian algorithm of \cite{PreSe:25} with a number of runs $R>1$, one is exposed to the risk of missing details in the coarse runs, causing a high error and a long time to converge. This is illustrated for $R=2$ in the right panel of Figure~\ref{fig:cartesianComparison} and can of course be avoided by employing finer grids from the beginning. Next, we have applied the Cartesian algorithm with $R=1$ and $\dx=\num{3.42e-2}$, mimicking the same changes of parameters and of reconstruction operators of the adaptive version. The central panel of Figure~\ref{fig:cartesianComparison} depicts the solution in this case and from Table~\ref{tab:cartesianComparison} one can appreciate that the error is comparable with the one of the novel adaptive approach, which is slightly faster; moreover, the Cartesian setup employs ten times more degrees of freedom thus requiring a much larger memory usage.

\subsection{Detecting cavities}
\label{ssec:tunnel}

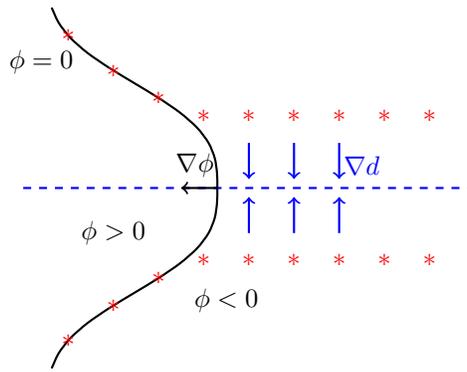
\begin{figure}[tpb]
\centering
\begin{tikzpicture}[scale=0.6]
    \coordinate (A) at (-5,0);
    \coordinate (B) at (5,0);
    \coordinate (C) at (0,0);

    \draw[blue, dashed, thick] (A) -- (B);

    \node[blue] at (2.5,0.5) {\(\nabla d\)};

    \node at (-0.5,-2.5) {\(\phi < 0\)};
    \node at (-4.6,2.8) {\(\phi = 0\)};
    \node at (-3,-1.) {\(\phi > 0\)};
    \node at (-1.2,0.5) {\(\nabla\phi\)};

    \draw[blue, thick, <-] (0,0.2) -- (0,1);
    \draw[blue, thick, <-] (1,0.2) -- (1,1);
    \draw[blue, thick, <-] (2,0.2) -- (2,1);

    \draw[blue, thick, <-] (0,-0.2) -- (0,-1);
    \draw[blue, thick, <-] (1,-0.2) -- (1,-1);
    \draw[blue, thick, <-] (2,-0.2) -- (2,-1);
    \draw[thick, ->] (-0.7,0) -- (-1.5,-0);

    \draw[domain=-4:4, smooth, variable=\x, black, thick] 
        plot ({3.8*exp(-(abs(\x*\x*\x))/19)-4.5},{\x});

    \node[red] at (-4,3.3777) {$*$};
    \node[red] at (-3,2.6042) {$*$};
    \node[red] at (-2,2.0) {$*$};
    \node[red] at (-1,1.6) {$*$};
    \node[red] at (-0,1.6) {$*$};
    \node[red] at (1,1.6) {$*$};
    \node[red] at (2,1.6) {$*$};
    \node[red] at (3,1.6) {$*$};
    \node[red] at (4,1.6) {$*$};
    \node[red] at (-4,-3.3777) {$*$};
    \node[red] at (-3,-2.6042) {$*$};
    \node[red] at (-2,-2.0) {$*$};
    \node[red] at (-1,-1.6) {$*$};
    \node[red] at (0,-1.6) {$*$};
    \node[red] at (1,-1.6) {$*$};
    \node[red] at (2,-1.6) {$*$};
    \node[red] at (3,-1.6) {$*$};
    \node[red] at (4,-1.6) {$*$};
\end{tikzpicture}
\vspace{0.3cm}
\caption{Configuration of a tunnel in a point cloud $\pcloud$, represented by the red star points. The evolution of the $\phi$ initially pushes its zero level set inside of the corridor, but it is then slowed down until it stops due to the vanishing of the transport term. In fact, once inside the tunnel, the gradient of $\phi$ and the velocity field $\nabla d$ become orthogonal.} 
\label{fig:tunnel}
\end{figure}

\begin{figure*}
\vspace{-1cm}
\centering
\begin{minipage}{0.45\linewidth}
    \centering    \includegraphics[width=1.\linewidth]{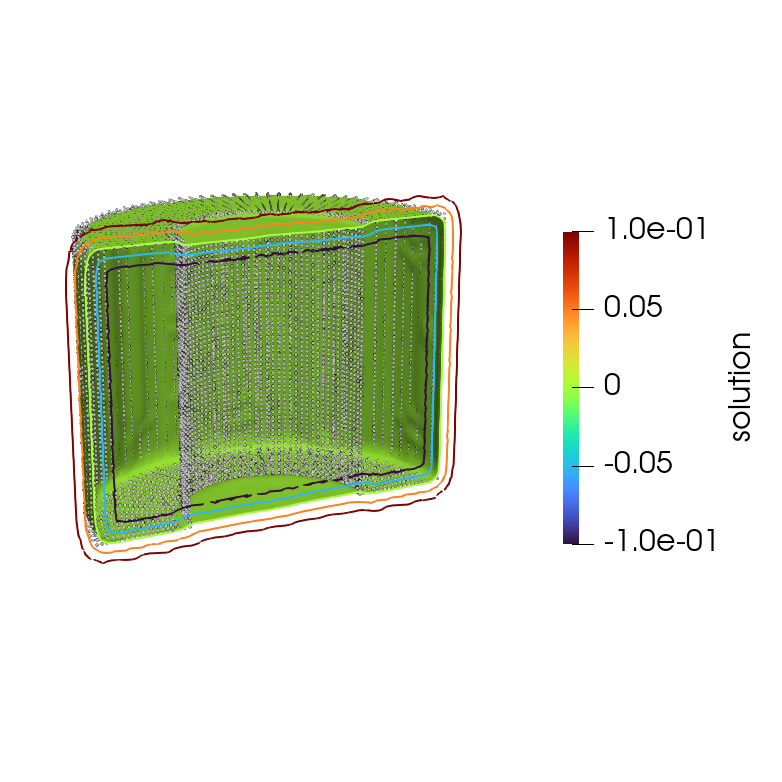}
\end{minipage}\hspace{0.2cm}
\begin{minipage}{0.45\linewidth}
    \centering
    \includegraphics[width=1.\linewidth]{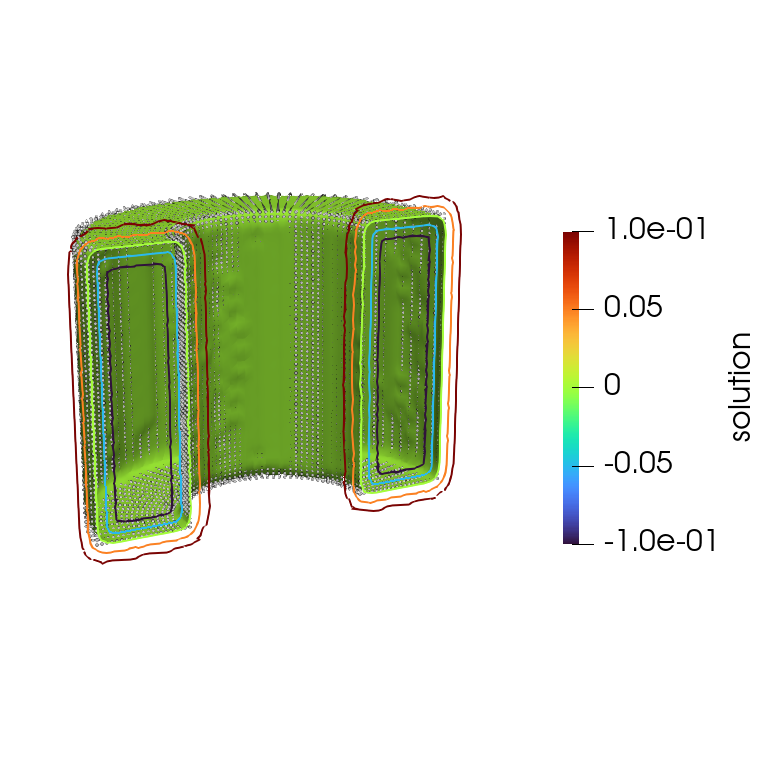}
\end{minipage}
\vspace{-1.2cm}
    \caption{From left to right: evenly spaced isocontours computed for the ``Cylinder'' test at the end of the evolution, using the original velocity and the modified one, respectively. Modifying the velocity field as in \S~\ref{ssec:tunnel}, the level set manages to detect the cavity.}
    \label{fig:cylinderIsocontours}
\end{figure*}

The last part of the numerical tests is devoted to particular shapes that present tunnel-type features corresponding to real cavities, which we would like to distinguish from fake ones due to missing data. A 2D example is depicted in Figure~\ref{fig:tunnel}, where one can notice that if a portion of the cloud $\pcloud$ is distributed along two close parallel lines, the evolution would be stopped by the vanishing of the transport term, namely by the orthogonality between $\nabla\phi$ and $\nabla d$, even if the distance function is still large. This is topologically analogous to the case of holes in the ``Knot'' test, but metrically different, since by ``tunnels'' we mean ``holes'' that are sufficiently long with respect to the size of their cross section. 

To address this issue we propose to modify the governing PDE \eqref{eq:levelset:pde} if the following conditions hold: on an active quadrant $j\in\tilde\grid$, if $d_j>4 h_\pcloud$ and $\modulo{\nabla d_j}<0.9$, we replace the velocity field $\nabla d$ with $\nabla\phi$, which is in practice a way of letting the level set $\phi$ proceed in the same direction followed so far. Note that a condition on the scalar product $\nabla\phi\cdot\nabla d$ might be too strong and might expose to the risk of perturbing the evolution even in well resolved areas.

We first consider a synthetic object constituted by a cylinder with a hole connecting the two parallel faces. The cloud has $16000$ points, its resolution is $h_\pcloud\approx\num{2.40e-2}$ and we start with $\dxmin=\num{3.19e-2}$, by setting $C_\pcloud=1$. The results obtained with the original PDE \eqref{eq:levelset:pde} and with the modified one are depicted in Figure~\ref{fig:cylinderIsocontours}, clearly showing the difference between the two algorithms in properly detecting the cavity at the centre.

\begin{figure*}
\vspace{-0.5cm}
\begin{minipage}{0.33\linewidth}
\includegraphics[width=1.3\linewidth]{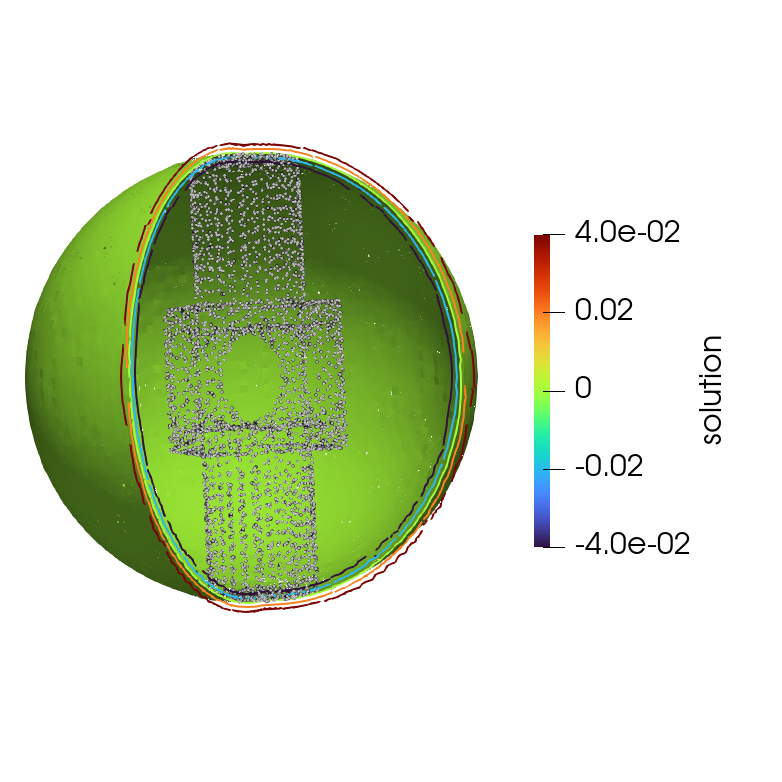}
\end{minipage}
\begin{minipage}{0.33\linewidth}
\includegraphics[width=1.3\linewidth]{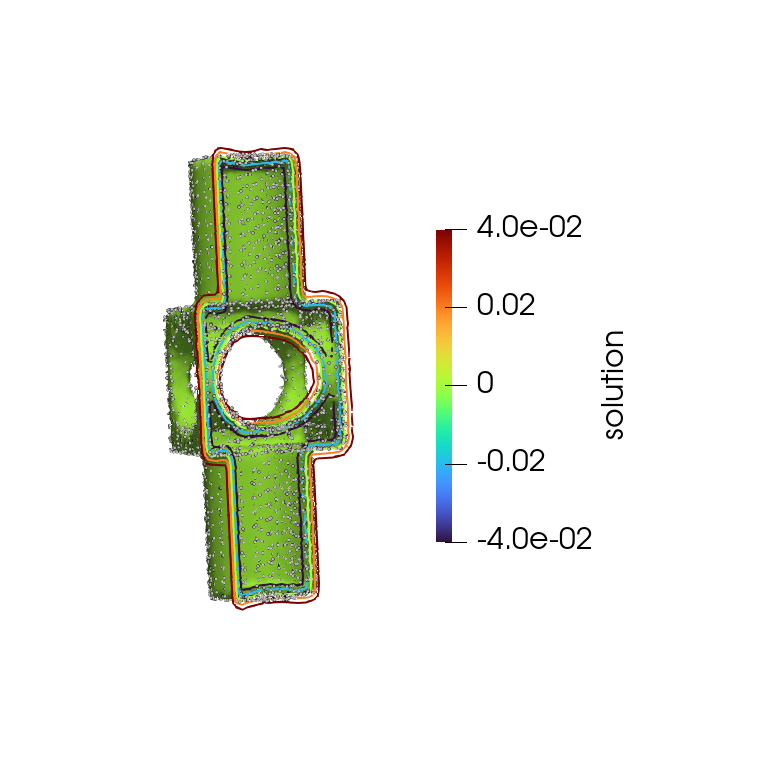}
\end{minipage}\hspace*{-0.8cm}
\begin{minipage}{0.33\linewidth} \includegraphics[width=1.2\linewidth]{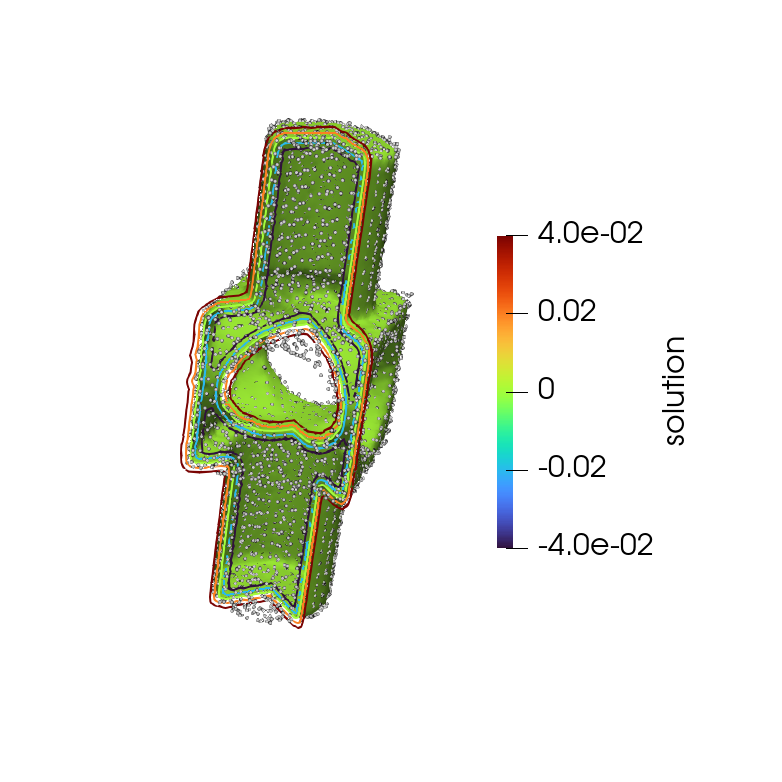}
\end{minipage}
\vspace{-0.6cm}
    \caption{From left to right: evenly spaced isocontours computed for the ``Mechaninal Part'' test at beginning and from two different perspectives at the end of the evolution. For this test, the modified velocity field, as described in \S~\ref{ssec:tunnel}, is considered and it allows the level set to detect the cavity in the centre.}
    \label{fig:mechPartReconstruction}
\end{figure*}

Similar considerations can be done for the ``Mechanical Part'' test: we process a point cloud of $4102$ points, with $h_\pcloud\approx\num{1.84e-2}$, presenting various features like sharp edges and corners, rounded parts, and a cavity in the centre. We set $C_\pcloud=0.5$ and thus set $\dxmin=\num{1.04e-2}$ allowing the modified velocity to come into play during the evolution. The result of the reconstruction procedure is shown in Figure~\ref{fig:mechPartReconstruction}.

\paragraph{Remark} 
We point out that the modification of the governing PDE proposed in this section can be applied for all the test cases presented in this work, even in presence of missing data, which is the case of the ``Frog'' of \S~\ref{ssec:realData}. Indeed, in such cases, the condition $\modulo{\nabla d_j}<0.9$ would not be satisfied and the computed level set will not protrude inside the object.

\section{Conclusions}
\label{sec:conclusions}

We presented a fully adaptive method for the reconstruction of surfaces from a set of unorganized points, having no information about their connection nor orientation. The LSM is used to detect and evolve these surfaces in an implicit way. In particular, we aim to keep the level set function close to the signed distance function, at least in the vicinity of its zero isocontour, in order to preserve numerical accuracy and to allow end users to reliably compute surface normals and curvature from the level set itself. The numerical method is based on a SL scheme on a quadtree (octree in 3D) grid and is coupled with a $\Pone$ or $\CWENO$ reconstruction operator based on the least-squares approach.
We made use of the grid and parallel implementation of the P4EST library, and presented all the complementary procedures needed to design the algorithm.

Future research will be directed toward a finer and local estimation of the parameters employed in the method, taking into account the curvature of the level set function $\phi$ and a local estimation of the cloud size $h_{\pcloud}$, aiming to save computational effort, without compromising the quality of the results. Also, employing a forest of trees, instead of a single octree, would be more efficient.

\paragraph{Acknowledgements}

The authors acknowledge the CINECA award under the ISCRA initiative, for the availability of high-performance computing resources and support (ISCRA Projects ID: HP10C7HWOL, HP10CO8NC7, HP10COTGT4). Both authors are members of the Gruppo Nazionale Calcolo Scientifico-Istituto Nazionale di Alta Matematica (GNCS-INdAM). This research has been partly funded by the PRIN-PNRR project ``MATHematical tools for predictive maintenance and PROtection of CULTtural heritage (MATHPROCULT)'' (code P20228HZWR).

\paragraph{Author Contributions} Both authors contributed equally to the conceptualization of the method. Silvia Preda has been the primary implementer of the method, designed and performed the numerical tests.

\paragraph{Data Availability} The datasets analysed during the current study are available in the repositories cited in the text.

\section*{Declarations}

The authors declare that they have no conflict of interest.

$\,$

$\,$

\bibliographystyle{spmpsci}
\bibliography{levelset.bib}
\end{document}